\documentclass[11pt]{article}

\usepackage{amsmath, amsfonts, amssymb,latexsym,wasysym,graphicx,stmaryrd,tikz,txfonts}
\input epsf.sty

\newtheorem{Theorem}{Theorem}[section]
\newtheorem{Lemma}{Lemma}[section]
\newtheorem{Proposition}[Lemma]{Proposition}
\newtheorem{Corollary}[Lemma]{Corollary}
\newtheorem{Definition}[Lemma]{Definition}

\newtheorem{Assumption}[Lemma]{Assumption}

\newcommand{\BEQ}{\begin{equation}}     
\newcommand{\BEA}{\begin{eqnarray}}
\newcommand{\BD}{\begin{displaymath}}
\newcommand{\EEQ}{\end{equation}}       
\newcommand{\EEA}{\end{eqnarray}}
\newcommand{\ED}{\end{displaymath}}

\newcommand{\del}{\delta}
\newcommand{\Del}{\Delta}
\newcommand{\eps}{\varepsilon}          

\newcommand{\supp}{{\mathrm{supp}}}

\newcommand{\Vol}{{\mathrm{Vol}}}

\newcommand{\R}{\mathbb{R}}

\newcommand{\Z}{\mathbb{Z}}
\newcommand{\N}{\mathbb{N}}
\newcommand{\D}{\mathbb{D}}

\def\proba{{\mathbb{P}}}

\def\esper{{\mathbb{E}}}

\def\locsup{{\mathrm{loc}}\, {\mathrm{sup}}}
\def\Locsup{{\mathrm{Loc}}\, {\mathrm{sup}}}

\newcommand{\eop}{\hfill $\Box$}        

\newcommand{\half}{{1\over 2}}          
\renewcommand{\vec}[1]{\boldsymbol{#1}} 

\catcode`\@=11
\def\numberbysection{\@addtoreset{equation}{section}
        \def\theequation{\thesection.\arabic{equation}}}
\numberbysection

\begin{document}

\vspace*{1.5cm}
\begin{center}
{\Large \bf  Generalized PDE estimates for KPZ equations through 
Hamilton-Jacobi-Bellman formalism}

\end{center}

\vspace{2mm}
\begin{center}
{\bf  J\'er\'emie Unterberger$^a$}
\end{center}

\vskip 0.5 cm
\centerline {$^a$Institut Elie Cartan,\footnote{Laboratoire 
associ\'e au CNRS UMR 7502} Universit\'e de Lorraine,} 
\centerline{ B.P. 239, 
F -- 54506 Vand{\oe}uvre-l\`es-Nancy Cedex, France}
\centerline{jeremie.unterberger@univ-lorraine.fr}

\vspace{2mm}
\begin{quote}

\renewcommand{\baselinestretch}{1.0}
\footnotesize
{We study in this series of articles the Kardar-Parisi-Zhang  (KPZ) equation
$$ \partial_t h(t,x)=\nu\Del h(t,x)+\lambda  V(|\nabla h(t,x)|) +\sqrt{D}\, \eta(t,x), \qquad x\in\R^d $$
in $d\ge 1$ dimensions. The forcing term $\eta$ in the right-hand side is
 a regularized white noise. The deposition rate $V$ is assumed to be isotropic and convex.
Assuming $V(0)\ge 0$, one finds $V(|\nabla h|)\ltimes |\nabla h|^2$ for small gradients,
 yielding the equation which is  most commonly used  in the literature. 
 
\medskip
The present article, a continuation of \cite{Unt-KPZ1},  is dedicated to a generalization of the  PDE estimates obtained in the previous
article  to the case of a deposition rate $V$ with polynomial growth of arbitrary order at infinity, for which in general the Cole-Hopf transformation does not allow any more
a comparison to the heat equation. The main tool here instead
is the representation of $h$ as the solution of some minimization problem through the Hamilton-Jacobi-Bellman formalism. This sole representation turns out to be powerful enough to produce {\em local} or {\em pointwise}  estimates in ${\cal W}$-spaces
of functions with "locally bounded averages", as in \cite{Unt-KPZ1}, implying in 
particular global existence and uniqueness of solutions. 
}
\end{quote}

\vspace{4mm}
\noindent


\smallskip


\newpage

\tableofcontents



\section{Introduction and statement of main results}



\subsection{Introduction}


We consider in this series of articles inhomogeneous, non-linear, viscous Hamilton-Jacobi equations of the form
\BEQ \partial_t h(t,x)=\Del h(t,x)-\eps h(t,x)+\lambda V(\nabla h(t,x))+g(t,x), \qquad h(0,x)=h_0(x) \label{eq:KPZ} \EEQ
on $\R_+\times\R^d$, $d\ge 1$, 
where $h_0$ is the initial condition, $-\eps h$ $(\eps\ge 0)$   a linear damping term, $g$ a forcing term, and $\lambda>0$ a positive constant. The non-negative function  $V:\R^d\to\R_+$, called
{\em deposition rate}, is only  assumed here to be convex, besides  very general properties ($C^2$ regularity, 
polynomial growth at infinity). Physically, the above equation modelizes the growth of an interface under (i) diffusion; (ii)
material deposition at site $x$ depending only on the gradient of the interface at that point; (iii) a forcing term
$g$ viewed as a noise.  By  a vertical drift,
$h(t,x)\to h(t,x)-tV(0)$, and a change of coordinate, $x\to x-t\nabla V(0)$, one may (and shall) assume
that $V(0)=0$ and $\nabla V(0)=0$, so that $V(\nabla h)=O(|\nabla h|^2)$ at a site where the interface is locally
almost flat (horizontal), i.e. for $|\nabla h|$ small. For physical reasons (although this condition is by no means necessary
for the estimates developed in this article), we shall also assume $V$ to be isotropic. 

Such PDEs are generalizations of the {\em Kardar-Parisi-Zhang} (KPZ) equation,
\BEQ \partial_t h(t,x)=\Del h(t,x)+\lambda \left(\sqrt{1+|\nabla h(t,x)|^2} - 1\right)+g(t,x), \qquad h(0,x)=h_0(x)\EEQ
also written (using second-order Taylor expansion around a locally flat interface),
\BEQ \partial_t h(t,x)=\Del h(t,x)+\frac{\lambda}{2} |\nabla h(t,x)|^2+g(t,x), \qquad h(0,x)=h_0(x).\EEQ

As explained in details in the introduction of the previous article \cite{Unt-KPZ1}, (1) {\em PDE estimates for multi-dimensional KPZ equation},  our general motivation is to study, by a rigorous
implementation of K. Wilsons's renormalization scheme, the large-scale limit in dimensions $d\ge 3$ of the {\em noisy}
KPZ equation, for which $g=\eta$ is a regularized white noise. We are  however only concerned in {\em single-scale} estimates here, so the dimension plays no r\^ole for
the moment, and our results apply equally well in case $d=1,2$.
 The {\em scale} $j$ {\em infra-red cut-off KPZ equation} $(j\ge 0)$ is given by (\ref{eq:KPZ}) with $\eps=2^{-j}$, and
 a right-hand side $g$ satisfying bounds typical for time-derivatives of 'averaged functions' of the form $e^{t\Del}f$ with $t\ge 2^j$, roughly speaking,
 $|\nabla^{\kappa} g|=O((2^{-j})^{1+|\kappa|/2})$, as follows from standard parabolic estimates.  If $h_0$, $\nabla h_0$, $g$  and $\nabla g$ are {\em bounded}, then the maximum and comparison principles apply to solutions of eq. (\ref{eq:KPZ}). As a matter of fact,
 a lot is known ($L^p$-bounds for $h$ or its gradient, asymptotic long-time behaviour,...) about the solutions
 of the homogeneous equation ($g=0$), see e.g. \cite{BenBenLau}, \cite{BenKarLau}, \cite{BenLau}, \cite{GilGueKer},
 \cite{LauSou}; although solutions are smooth for $t>0$, the theory of viscosity solutions plays an important r\^ole
 in these developments. 
 
We are however typically interested in {\em unbounded} solutions which arise naturally when $g$ is a
{\em space-translation invariant in law}, random forcing term, and thus has {\em large,
unbounded fluctuations}. With the application to the noisy KPZ equation in
mind, we introduced in \cite{Unt-KPZ1} new functional spaces ${\cal W}_j^{1,\infty;\lambda}\subset W^{1,\infty}_{loc}(\R^d)$ and their time-dependent versions,
${\cal W}_j^{1,\infty;\lambda}([0,T])\subset C([0,T],W^{1,\infty}_{loc}(\R^d))$ ($T>0$) having roughly the following properties (see below "Notations" if needed):

\begin{itemize}
\item[(i)] ${\cal W}^{1,\infty}\subset {\cal W}_j^{1,\infty;\lambda}$, $C([0,T],{\cal W}^{1,\infty})\subset {\cal W}^{1,\infty;\lambda}_j([0,T])$, where ${\cal W}^{1,\infty}\subset L^{\infty}(\R^d)$ is the subspace of functions
with bounded generalized gradient. 
\item[(ii)] ${\cal W}_j^{1,\infty;\lambda}$ (and
similarly ${\cal W}_j^{1,\infty;\lambda}([0,T])$) is given in terms of a family of {\em local quasi-norms} $|||\ \cdot \ |||_{{\cal W}_j^{1,\infty;\lambda}}(x)$, namely (see \cite{Unt-KPZ1}, \S 3.2 and 3.3), ${\cal W}_j^{1,\infty;\lambda}:=\{f\in W^{1,\infty}_{loc}(\R^d); \forall x\in\R^d,
|||f|||_{{\cal W}_j^{1,\infty;\lambda}}(x)<\infty\}$, and 
\BEQ |||f|||_{{\cal W}_j^{1,\infty;\lambda}}(x)\le |||f|||_{{\cal W}_j^{1,\infty;\lambda'}}(x) \qquad (\lambda\le\lambda'), \qquad
|||\mu f|||_{{\cal W}_j^{1,\infty;\lambda}}(x)\le |\mu| \  |||f|||_{{\cal W}_j^{1,\infty;|\mu|\lambda}}(x) \qquad (\mu\in\R);\EEQ
\BEQ |||f_1+f_2|||_{{\cal W}_j^{1,\infty;\lambda}}(x)\le \frac{1}{p_1}|||p_1 f_1|||_{{\cal W}_j^{1,\infty;\lambda}}(x)+\frac{1}{p_2}|||p_2 f_2|||_{{\cal W}_j^{1,\infty;
\lambda}}(x) \qquad (p_1,p_2\ge 1,\ \frac{1}{p_1}+\frac{1}{p_2}=1),\EEQ
from which it follows in particular that ${\cal W}_j^{1,\infty;\lambda}$ is a convex subset of $W^{1,\infty}_{loc}(\R^d)$.
 \item[(iii)] A comparison principle holds
for solutions in ${\cal W}_j^{1,\infty;\lambda}$, at least for the homogeneous
KPZ equation (obtained by letting $g=0$). 
\item[(iv)] Let $\lambda'>\lambda$. Then eq. (\ref{eq:KPZ}) with $\eps=2^{-j}$ has a unique classical solution
$h$ such that $h_t\in {\cal W}_j^{1,\infty;2\lambda/5}$ for all $t\in[0,T]$ if $h_0\in {\cal W}_j^{1,\infty;2\lambda'}\cap C^2(\R^d)$ and
$g\in {\cal W}_j^{1,\infty;2\lambda}([0,T])\cap C([0,T],C^3(\R^d))$. One has an explicit, $t$-independent bound on $|||h_t|||_{{\cal W}_j^{1,\infty;2\lambda/5}}(x)$ in terms of $|||h_0|||_{{\cal W}_j^{1,\infty;2\lambda'}}(x)$ and $|||g|||_{{\cal W}_j^{1,\infty;2\lambda}([0,t])}(x)$. Briefly said, if the data $h_0,g$ are sufficiently smooth,
then one has a {\em local bound} on $h_t(x)$ in local ${\cal W}$-quasi norm in terms of the
relevant local ${\cal W}$-quasi norms of $h_0$ and $g$ {\em at the same point}, $x$. 
\end{itemize}
The general principles underlying the definition of these functional spaces are retraced in \S 0.3 below.

\medskip

\noindent One of the main ingredients in the proof of these results has been the {\em Cole-Hopf} (or {\em exponential}) {\em transformation}, $h\mapsto e^{\lambda h}$.
This transformation maps a solution of (\ref{eq:KPZ}) into a sub-solution of the linear heat equation provided 
$V(|\nabla h|)\le |\nabla h|^2$.  The above 
pointwise quasi-norms measure {\em local averages of the Cole-Hopf transform} of their
arguments. If $V$ is not quadratically bounded at infinity, essentially all our conclusions
in \cite{Unt-KPZ1} fail.

\bigskip
\noindent We tackle here the same questions from a different perspective, starting from a {\em Hamilton-Jacobi-Bellman}
representation of the solutions of (\ref{eq:KPZ}): the classical idea (see section 1 for a thorough discussion) is to rewrite
the function $h_t$  as the maximum over an admissible class of random paths $X$ driven by Brownian motion $B$ of
a functional $\int_0^t F(s,X_s)ds$.  To be explicit, letting $\tilde{V}$ be the 
Legendre transform of $V$, we have (see Lemma \ref{lem:HJB-representation})
$ h(t,x)=\sup_{\alpha} J_{\alpha}(t,x)$,  where
\BEQ J_{\alpha}(t,x):=\esper^{0,x} \left[  \int_0^t e^{-2^{-j}s} \left(-\lambda \tilde{V}(\frac{\alpha_s}{\lambda}) 
+g(t-s,X^{\alpha}_s) \right) ds + e^{-2^{-j} t} h_0(X^{\alpha}_t)\right],  \EEQ
 $(X^{\alpha}_s)_{0\le s\le t} $ being the solution of the stochastic differential equation, $dX^{\alpha}_s=\alpha_s ds+ dB_s$ with
initial condition $X^{\alpha}_0=x$.

\medskip
\noindent Assuming $V$ is polynomially bounded at infinity (see \S 0.2 below for
precise assumptions), $\tilde{V}(\alpha_s)$ is bounded below by some positive power of 
$\alpha_s$. Note also that $X^{\alpha}$ is essentially bounded by the integral of $\alpha$,
except if $B$ is larger. Thus we have a dichotomy: (i) either the integral of $\alpha$ is large, and $J_{\alpha}(t,x)$ is small; (ii) or the integral of $\alpha$ is negligible
with respect to $B$, and $J_{\alpha}(t,x)\simeq J_0(t,x)$ is essentially given by the solution of the linear heat equation $(\partial_t-\Del+2^{-j})h=g$.
This random control-theoretic argument can be made quantitative, and yields a priori estimates for the
KPZ solution (see Theorem \ref{th:main}).

\bigskip
\noindent We postpone a detailed discussion of our main results (see \S 0.4) and
a brief outline of the article (see \S 0.5). The reason is two-fold: first, though
techniques differ widely, and our hypotheses are somewhat different also,  the 
general philosophy, and our final statements, look pretty much the same as in our previous article,
 to which we therefore refer the reader. Second, there are some relatively minor, but technically important,
 differences between the time-dependent spaces of \cite{Unt-KPZ1} and those of the
 present work; also, our techniques make it possible to work with more general (and
 larger) ${\cal W}$-spaces than the Cole-Hopf ${\cal W}$-spaces mentioned above.
 So one definitely needs a precise definition of the functional spaces (see \S 0.3) to
 appreciate the originality of this paper with respect to the previous one.
 
\medskip
\noindent The reader may be interested in comparing these results with those we
 obtained in \cite{Unt-Bur2} for unbounded solutions of Burgers' equation. The estimates
 we get in that case are also based on a random characteristic representation, but the
 way fluctuations are handled is very dissimilar.

\medskip
{\bf Notations.}  The notation:
$f(u)\lesssim g(u)$, resp. $f(u)\gtrsim g(u)$ means: $|f(u)|\le C|g(u)|$, resp. $|f(u)|\ge C|g(u)|$, where $C>0$ is an unessential constant. 
Similarly, $f(u)\approx g(u)$ means: $f(u)\lesssim g(u)$ and $g(u)\lesssim f(u)$.
We denote by $L^{\infty}$, resp. $L^{\infty}_{loc}$  the space of bounded, resp.
locally bounded functions on $\R^d$ or $\R_+\times\R^d$;  by ${\cal W}^{1,\infty}$, resp. ${\cal W}^{1,\infty}_{loc}$  the  space of bounded, resp. locally bounded  functions on $\R^d$ or $\R_+\times\R^d$  with bounded, resp. locally generalized derivative;  and by $C^{1,2}$ the space of functions which are $C^1$ in time and $C^2$ in space. The {\em average} 
$\frac{1}{|\Omega|} \int_{\Omega} f$ of a function $f$
on a bounded domain $\Omega\subset\R^d$ is denoted by $\fint_{\Omega}f$.


\subsection{The model}


A {\em KPZ equation}  is a viscous Hamilton-Jacobi
equation,
\BEQ \partial_t h(t,x)=\Del h(t,x) +\lambda V(\nabla h(t,x))+g(t,x), \label{eq:full-KPZ} \EEQ
where $g$ is a right-hand side (or {\em forcing term}) in a suitable functional space. 
For simplicity we have chosen the diffusion constant $\nu$ to be $1$ (which can be done by a simple 
scaling).  The deposition rate $V$
satisfies in this article the following assumptions.

\medskip

\begin{Assumption} \label{assumptions} 
The deposition rate $V$ satisfies the following assumptions,
\begin{itemize}
\item[(1)] $V$ is $C^2$;
\item[(2)] $V$ is isotropic, i.e. $V(\nabla h)$ is a function of $y=|\nabla h|$; by abuse of notation we shall
consider $V$ either as a function of $\nabla h$ or of $y$;
\item[(3)] $V$ is convex;
\item[(4)] $V(0)=0$ and $0\le V(y)\le \max(\half y^2,\half y^{\beta})$ for all $y\ge 0$, where $\beta\ge 2$ is called the
{\em growth exponent at infinity} of $V$.
\end{itemize} 
\end{Assumption}

\medskip

Compared to the above class of equations, the {\em infra-red cut-off KPZ equation}, that
we shall  be dealing with exclusively in this paper, has a supplementary linear
damping term.  

\begin{Definition}[infra-red cut-off equation]
Let $V$ satisfy the above assumptions, and $j\in\N$. Then the {\em infra-red cut-off KPZ equation of scale $j$} with {\em data} $(h_0,g)$  is
the following  equation,
\BEQ \partial_t h(t,x)=\Del h(t,x) - 2^{-j}h(t,x) +\lambda V(\nabla h(t,x))+g(t,x),
\qquad h\big|_{t=0}=h_0. \label{eq:KPZdowntoj} \EEQ
\end{Definition}

The extra term $-2^{-j}h$ in the right-hand side implies in principle an exponential
decay of memory. As discussed in \cite{Unt-KPZ1}, sections 4 and 5, the operator $(\Del-2^{-j})^{-1}$ is a kind of ersatz for the 
high-momentum propagator $G^{j\to}:=\sum_{i=0}^j G^i$ with scale $j$ infra-red cut-off.


\subsection{Functional spaces}

The purpose of this article is to show that the infra-red cut-off  KPZ equation of scale $j$, 
eq. (\ref{eq:KPZdowntoj}), has
a single solution $h$ in an adequate functional space under suitable assumptions on the right-hand side $g$, and to
give appropriate bounds in suitable norms for $h_t$ in terms of $h_0$ and $g$. 

For the applications we have in mind (including the noisy KPZ equation, or more generally viscous Hamilton-Jacobi
equations with an extra noise term in the right-hand side), the initial condition $h_0$ and the right-hand side $g$
are unbounded, which led us to introduce new functional spaces, ${\cal W}^{1,\infty;\lambda}_j$
and ${\cal W}^{1,\infty;\lambda}_j([0,T])$ in the previous article. Compared to \cite{Unt-KPZ1}, our arguments here
have larger range of validity but give less precise bounds, which requires some rather minor changes in our
definitions. However, the general principles  underlying the construction of all these spaces is the same.
We now describe them briefly and refer to \cite{Unt-KPZ1}, sections 3 and 4 for
more details.

 Assume first $h_0\in{\cal W}^{1,\infty}$ and 
 $g\in L^1_{loc}(\R_+,{\cal W}^{1,\infty})$. By classical arguments derived from the parabolic maximum principle
 (see \cite{Unt-KPZ1}, section 2), the associated Cauchy problem (\ref{eq:full-KPZ}) has
a unique, global solution $h$ which lies in ${\cal W}^{1,\infty}$ for all $t\ge 0$ and  is classical for
strictly  positive times, that is, $h\in C([0,+\infty)\times\R^d)\cap
C^{1,2}((0,\infty)\times\R^d)$. Furthermore,
\BEQ ||h_t||_{\infty}\le ||h_0||_{\infty}+\int_0^t ||g_s||_{\infty} ds. \label{eq:maximum-principle} \EEQ

Recall however that the emphasis in this series of articles is in the {\em noisy} KPZ equation, for which $g=\eta$
 is a (suitably regularized) white noise.  Generally speaking,
$h$ is expected to behave more or less like the solution $\phi$ of the linearized equation,
\BEQ \partial_t \phi=\Del \phi +g, \label{eq:OU} \EEQ
which is (letting $h_0\equiv 0$ just for the sake of discussion) simply $
(\partial_t-\Del)^{-1}g=:Gg$. 
If $g=\eta$, then $g$ and $\phi=Gg$ will be space-translation invariant in law, hence are a.s. {\em unbounded}.
Thus we cannot expect $h_t$ to lie in ${\cal W}^{1,\infty}$. However, {\em local averages} 
\BEQ \fint_{B(x,r)} dy |\eta_t(y)|:= \frac{\int_{B(x,r)} dy |\eta_t(y)|}{\Vol(B(x,r))}, \qquad r>0 \EEQ
  of $\eta$ are {\em locally
uniformly bounded}, which amounts to saying that $\eta\in {\cal H}^0(\R^d)$, in the sense of the following definition.

\begin{Definition}
Let
\BEQ {\cal H}^0:=\{f\in L^{\infty}_{loc}(\R^d)\ |\ \forall x\in\R^d, f^*(x)<\infty\}\EEQ
where
\BEQ f^*(x):=\sup_{\tau>0} e^{\tau\Del}|f|(x)\in [0,+\infty].\EEQ
\end{Definition}

Due to the averaging and scaled decay properties  of the heat kernel, it is natural to expect that
$e^{\tau\Del}|f|(x)$ may be substituted with $\fint_{B(x,\sqrt{\tau})}dy |f(y)|$ in the above definition; it is actually proven in \cite{Unt-KPZ1}
that $f^*(x)\approx \sup_{r>0} \fint_{B(x,r)} dy |f(y)|$.

With this definition in hand, we may substitute the usual parabolic estimates,
\BEQ ||\nabla^k e^{t\Del}f||_{\infty}\lesssim t^{-\frac{k}{2}} ||f||_{\infty},\qquad k\ge 0 
\label{eq:parabolic-estimates} \EEQ

with the  stronger, {\em pointwise parabolic estimates} (see \cite{Unt-KPZ1}, section 3),
\BEQ (\nabla^k e^{t\Del}f)^*(x)\lesssim t^{-k/2}f^*(x), \qquad k\ge 0\EEQ
including the obvious but fundamental 
\BEQ (e^{t\Del}f)^*(x)\le f^*(x), \qquad t\ge 0.  \label{eq:obvious-but-fundamental} \EEQ

\medskip

Contrary to what was the case in the preceding article 
\cite{Unt-KPZ1} (see in particular Lemma 3.10), however, we cannot directly
compare KPZ solutions to solutions of the heat equation. The reason is that the
{\em Cole-Hopf transformation} $h\mapsto w:=e^{\lambda h}$ maps a solution of the
homogeneous KPZ equation $\partial_t h=\Del h+ \lambda V(|\nabla h|)$ to a solution
of the transformed KPZ equation $\partial_t w=\Del w+ \lambda (V(|\nabla h|)-|\nabla h|^2)w$. Now, {\em if} $V(|\nabla h|)\le |\nabla h|^2$ then $w$ is a sub-solution of the linear
heat equation. This property, which was one of the cornerstones in \cite{Unt-KPZ1}, fails
totally here. Thus  
(\ref{eq:obvious-but-fundamental})  cannot be used directly. 

As was already the case in \cite{Unt-KPZ1}, it is more convenient (although by no
means necessary) to cut-off small scale fluctuations of the initial condition
$h_0$, or of the right-hand side $g$,  by considering their {\em local suprema}.

\begin{Definition}[local supremum of order $j$] \label{def:locsup}
\begin{itemize}
 \item[(i)] Let $f:\R^d\to\R$ be a function in $L^{\infty}_{loc}(\R^d)$. Then $\locsup^j(f):\R^d\to\R$
is the function in $L^{\infty}_{loc}(\R^d)$ defined by
$\locsup^j(f)(x):=\sup_{y\in B(x,2^{j/2})} |f(y)|$.
\item[(ii)] (time-dependent case) Let $g:\R_+\times \R^d\to\R$ be a function in $L^{\infty}_{loc}(\R_+\times \R^d)$. Then $\Locsup^j(g):\R_+\times \R^d\to\R$
is the function in $L^{\infty}_{loc}(\R_+\times \R^d)$ defined by
$\Locsup^j(g)(t,x):=\sup_{s\in B(t,2^j)} \sup_{y\in B(x,2^{j/2})} |g(s,y)|$.
\end{itemize}
\end{Definition}

 This local supremum operation allows one to discretize space. Let
$\D^j$ be the set of all cells of the lattice $2^{j/2}\Z^d$, i.e. $\Del\in\D^j$ if and only
if $\Del=[2^{j/2}k_1,2^{j/2}(k_1+1)]\times \cdots\times [2^{j/2}k_d,2^{j/2}(k_d+1)]$ for
some $(k_1,\ldots,k_d)\in\Z^d$.  Clearly, if $r\gtrsim 1$,  
\BEQ r^{-d} \sum_{\Del\in \D^j;\ \Del\subset B(x,r2^{j/2})} \sup_{\Del}|f| \lesssim \fint_{B(x,r2^{j/2})} \locsup^j(f)(x)\lesssim   r^{-d} \sum_{\Del\in \D^j;\ \Del\cap B(x,r2^{j/2})\not=\emptyset} \sup_{\Del}|f|,\EEQ
so
\BEQ \locsup^j(f)^*(x):=(\locsup^j(f))^*(x) \approx \sup_n n^{-d} \sum_{\Del\in \D^j;\ \Del\subset B(x,n2^{j/2})} \sup_{\Del}|f| \EEQ
where $n$ ranges either over $\N^*$ or on the set of dyadic integers $2^k$, $k\in\N$. 
As already noted in section 3.3 of \cite{Unt-KPZ1}, 
\BEQ \sup_{\Del} |f|\lesssim d^j(0,\Del)^d \locsup^j(f)^*(0)  \label{eq:growth-condition} \EEQ
for $d^j(0,\Del):= 2^{-j/2}\min_{x\in\Del}|x|\gtrsim 1$, so a function $f$ such
that $\locsup^j(f)^*(0)<\infty$ has at most polynomial growth of order $d$ at infinity.
This is an explicit {\em limitation to the fluctuations} of $f$.

\medskip
Functions like $\locsup^j(f)^*(x)$ bound very efficiently expectations like
$\esper[|f(X)|]$ when the random variable $X$ (coming from the Hamilton-Jacobi-Bellman representation) has a known Gaussian queue: 

\begin{Lemma} \label{lem:integration0}
Let $f\in L^{\infty}_{loc}(\R^d)$ such that $\locsup^j(f)^*(x)<\infty$, and 
$X$ be an $\R^d$-valued random variable such that $|X-x|\le |B_t|$. Then
\BEQ \esper[|f(X)|] \lesssim (2^{-j}t)^{d/2} \locsup^j(f)^*(x). \EEQ
\end{Lemma}

{\bf Proof.} We bound $\esper[|f(X)|]$ by the sum over squares $\sum_{\Del\in\D^j}
\proba[X-x\in\Del] \sup_{\Del}|f|$, and  $\proba[X-x\in\Del]$ by
$\proba[|B_t|\ge 2^{j/2} d^j(0,\Del)]\le \sup_{y\in\Del} e^{- c|y|^2/t}$ for some $c>0$. Thus
$\esper[|f(X)|] \lesssim (2^{-j}t)^{d/2} e^{c't\Del}(\locsup^j(f))(x)$.   \hfill \eop

This key lemma, with its various generalizations, see Lemma \ref{lem:local-integration}  
and Lemma \ref{lem:pointwise-integration},  lies at the core of the proof of bounds of the KPZ solution in {\em local}
${\cal W}$-"quasi norms". Several families of these local "quasi-norms" have been introduced
in our previous article \cite{Unt-KPZ1}. It is useful at this point to recall how these
were defined. First
$|||\, \cdot \,|||_{{\cal H}^{\lambda}}(x)$ is essentially obtained from the "star" local quasi-norm
$f\mapsto f^*(x)$ through the conjugation by the Cole-Hopf exponential
transform $f\mapsto e^{\lambda f}$,
namely, $|||f|||_{{\cal H}^{\lambda}}(x):=\frac{1}{\lambda} \ln (e^{\lambda |f|})^*(x).$
The  analogue of these for a time-dependent function $g$  is
\BEQ |||g|||_{{\cal H}^{\lambda}([0,t])}(x):=2^{-j} \int_0^t ds\,  e^{-2^{-j}s} |||\ 2^j 
g(t-s,\cdot)\, |||_{{\cal H}^{\lambda}}(x). \label{intro:eq:g-Hlambda} \EEQ
Now we let:

\begin{itemize}
\item[(i)] {\em (local quasi-norm for the initial condition)}

 $|||h_0|||_{{\cal W}_j^{1,\infty;\lambda}}(x):= \max \left( |||\locsup^j h_0|||_{{\cal H}^{\lambda}}(x),
|||2^{j/2} \locsup^j |\nabla h_0|\ |||_{{\cal H}^{\lambda}}(x) \right)$;
\item[(ii)] {\em (time-dependent local quasi-norm for the right-hand side)}

$|||g|||_{{\cal W}_j^{1,\infty;\lambda}([0,t])}(x):=\max \left( |||\Locsup^j g|||_{{\cal H}^{\lambda}([0,t])}(x),
|||2^{j/2} \Locsup^j |\nabla g|\ |||_{{\cal H}^{\lambda}([0,t])}(x) \right)$.
\end{itemize}

The main result of \cite{Unt-KPZ1} (Theorem 2), an exact bound for 
$|||h_t|||_{{\cal W}_j^{1,\infty;\lambda}}(x)$,  may be cited in this weaker form:
there exists a constant $c\in(0,1)$ such that $|||h_t|||_{{\cal W}_j^{1,\infty;c\lambda}}(x)
\lesssim e^{-2^{-j}t} |||h_0|||_{{\cal W}_j^{1,\infty;\lambda}}(x)+ |||g|||_{{\cal W}_j^{1,\infty;\lambda}([0,t])}(x).$ Note the "loss of regularity" in the $\lambda$-exponent. 

In the present setting, the Cole-Hopf transform does not play any particular r\^ole, and
it is both natural and useful to enlarge the above family of functional spaces by allowing
conjugation by other functions than exponentials. Thus we define:

\begin{Definition}
Let, for $P,P_{\pm}:\R_+\to\R_+$ convex and strictly increasing with $P_-\le P_+$, and
$f:\R^d\to\R$, resp. $g:\R_+\times\R^d\to\R$,
\BEQ |||f|||_{{\cal H}^P}(x):= P^{-1} \left( (P(|f|))^*(x) \right);
\label{intro:def:f-norm} \EEQ
\BEQ |||g|||_{{\cal H}^{P_{\pm}}([0,t])}(x):=2^{-j}  P_+^{-1}\circ P_- \left(\int_0^t e^{-2^{-j}s}
 P_-^{-1}\circ P_+ \left(|||\ 2^j g(t-s,\cdot)\, |||_{{\cal H}^{P_+}}(x) \right) \right)
 \label{intro:def:g-norm} \EEQ
\end{Definition}

With this definition we obtain generalizations of the convexity properties (ii)
stated in \S 0.1, namely, considering e.g. only the family of
local quasi-norms $|||\, \cdot\, |||_{{\cal H}^P}(x)$, 
\BEQ |||f|||_{{\cal H}^{P^-}}(x)\le |||f|||_{{\cal H}^{P^+}}(x), \qquad P_-\le P_+; \EEQ
\BEQ |||\mu f|||_{{\cal H}^P}(x)\le P^{-1} \left( (1-\mu)P(0)+\mu (P(f))^*(x) \right)
\le |||f|||_{{\cal H}^P}(x); \EEQ
\BEQ ||| f_1+f_2 |||_{{\cal H}^P}(x)\le P^{-1} \left( \frac{1}{p_1} (P(p_1 f_1))^*(x)+ \frac{1}{p_2}
(P(p_2 f_2))^*(x) \right), \qquad p_1,p_2\ge 1, \frac{1}{p_1}+\frac{1}{p_2}=1
\EEQ
In particular, the functional spaces ${\cal H}^P:=\{f\in L^{\infty}_{loc}(\R^d) \ |\  \forall x\in \R^d,
|||f|||_{{\cal H}^P}(x)<\infty\}$ and ${\cal H}^{P_{\pm}}([0,t]):=\{f\in C([0,t],
L^{\infty}_{loc}(\R^d))\ |\ \forall x\in\R^d, |||g|||_{{\cal H}^{P_{\pm}}([0,t]]}(x)<\infty\}$ are {\em convex}.

\medskip
\noindent Our previous families of quasi-norms are thus particular examples of the more general
\BEQ   |||h_0|||_{{\cal W}_j^{1,\infty;P}}(x):= \max \left( |||\locsup^j h_0|||_{{\cal H}^{P}}(x),
|||2^{j/2} \locsup^j |\nabla h_0|\ |||_{{\cal H}^{P}}(x) \right), 
\label{eq:Wj1inftyP} \EEQ 
\BEQ \qquad\qquad
|||g|||_{{\cal W}_j^{1,\infty;P_{\pm}}([0,t])}(x):=\max \left( |||\Locsup^j g|||_{{\cal H}^{P_{\pm}}([0,t])}(x),
|||2^{j/2} \Locsup^j |\nabla g|\ |||_{{\cal H}^{P_{\pm}}([0,t])}(x) \right) 
\label{eq:Wj1inftyP0t} \EEQ
for $h_0:\R^d\to\R$, resp. $g:\R_+\times\R^d\to\R$  (compare with (i), (ii) above),
defining new functional spaces, ${\cal W}_j^{1,\infty;P}, {\cal W}_j^{1,\infty;P_{\pm}}([0,t])$.

\bigskip
\noindent  Now the {\em general idea} is to prove  bounds of the form
\BEQ |||h_t|||_{{\cal W}^{1,\infty;P_-}}(x) \lesssim e^{-2^{-j}t}
P_-^{-1}\circ P_+\left( |||h_0|||_{{\cal W}^{1,\infty;P_+}}(x) \right) + 
P_-^{-1}\circ P_+\left(|||g|||_{{\cal W}^{1,\infty;P_{\pm}}([0,t])}(x) \right)\EEQ
for some couples $(P_-,P_+)$ with $P_-\le P_+$, in which the "loss of regularity"
phenomenon is apparent. In our previous article we have obtained  bounds {\em in the exponential case}, i.e. for $P_-(z)=e^{c\lambda z}$, $P_+(z)=e^{\lambda z}$. Note that $P_-^{-1}\circ P_+$, $P_+^{-1}\circ P_-$ amount to a simple multiplication by a constant in this special case. Here we state new bounds (i) in the {\em exponential case}, but also (ii) in the {\em polynomial case} where $P_{\pm}(z)=
z^{d_{\pm}}$, with $d_-\ge 1$ and $d_+ - d_-$ large enough.  Note that now $P_-^{-1}\circ P_+(z)=z^{d_+/d_-}$, giving a "Besov" flavour to the space ${\cal H}^{P_{\pm}}([0,t])$
(see \cite{Trie}). On the other hand ${\cal H}^P$, ${\cal H}^{P_{\pm}}([0,t])$ are
easily shown to be {\em linear subspaces} in this case, contrary to case (i).


\subsection{Main results}


For technical reasons, the time-dependent ${\cal W}$-spaces must be slightly modified 
(see comment after Theorem 2).
For $d'\ge 0$, we let
\BEQ |||g|||_{\widetilde{\cal W}_j^{1,\infty;P_{\pm}}([0,t])}(x):=\max \left( |||\Locsup^j g|||_{\widetilde{\cal H}^{P_{\pm}}([0,t])}(x),
|||2^{j/2} \Locsup^j |\nabla g|\ |||_{\widetilde{\cal H}^{P_{\pm}}([0,t])}(x) \right), \EEQ
with 
\BEQ |||g|||_{\widetilde{\cal H}^{P_{\pm}}([0,t])}(x):=2^{-j}  P_+^{-1}\circ P_- \left(\int_0^t e^{-2^{-j}s} (1+2^{-j}s)^{d'/2}
 P_-^{-1}\circ P_+ \left(|||\ 2^j g(t-s,\cdot)\, |||_{{\cal H}^{P_{\pm}}}(x) \right) \right)
\label{intro:tilde-g-norm} \EEQ 
(compare with (\ref{intro:def:g-norm}),(\ref{eq:Wj1inftyP0t})), defining as above new
functional spaces, $\widetilde{\cal H}^{P_{\pm}}([0,t])$ and $\widetilde{\cal W}_j^{1,\infty;P_{\pm}}([0,t])$. We leave intentionnally the dependence on the exponent $d'$
implicit.

\medskip\noindent
Let us first state the

\begin{Definition}[${\cal W}_j^{1,\infty;P_-}$-solution] \label{intro:def:Wj1inftyPsolution}
Choose $P_{\pm}:\R_+\to\R_+$ be two convex, strictly increasing functions with $P_-\le P_+$. 
Let $h_0\in {\cal W}_j^{1,\infty;P_{+}}$, $g\in \widetilde{\cal W}_j^{1,\infty;P_{\pm}}([0,T])$ and $h\in C([0,T];\widetilde{\cal W}_j^{1,\infty;P_-})$. The function $h$
is said to be a {\em ${\cal W}_j^{1,\infty;P_-}$-solution} of the scale $j$ infra-red cut-off KPZ equation with
right-hand side $g$ and initial condition $h_0$ if there exists a sequence of functions $h_0^{(n)}\in {\cal W}^{1,\infty}$, 
$g^{(n)}\in L^{\infty}([0,T];{\cal W}^{1,\infty})$ such that (i) for every compact $K\subset\R^d$, 
$h_0^{(n)}\to_{n\to\infty} h_0$ in ${\cal W}^{1,\infty}(K)$ and 
$g^{(n)}\to_{n\to\infty} g$ in $C([0,T];{\cal W}^{1,\infty}(K))$;
(ii) there exists a constant $C>1$ such that, for all $n\ge 1$ and $x\in\R^d$,
$|||h_0^{(n)}|||_{{\cal W}_j^{1,\infty;P_+}}(x)\le C |||h_0|||_{{\cal W}_j^{1,\infty;P_+}}(x)$,
 $|||g^{(n)}|||_{\widetilde{\cal W}_j^{1,\infty;P_{\pm}}([0,t])}(x)\le 
C |||g|||_{\widetilde{\cal W}_j^{1,\infty;P_{\pm}}([0,t])}(x)$;  (iii) $h^{(n)}\to_{n\to\infty} h$ in
$C([0,T];{\cal W}^{1,\infty}(K))$, where $h^{(n)}$ is the unique classical solution in $C([0,T];{\cal W}^{1,\infty})\cap
C^{1,2}((0,T]\times\R^d)$ of the KPZ equation 
\BEQ ({\mathrm{KPZ}}_n): (\partial_t-\Del) h^{(n)}=-2^{-j} h^{(n)}+\lambda V(\nabla h^{(n)})+g^{(n)},\qquad h^{(n)}(t=0)=h_0^{(n)}.\EEQ

\end{Definition}

\bigskip
\noindent
Our main results are the following.

\bigskip
\noindent
{\bf Theorem 1 (existence and unicity).} 
{\em Let $h_0\in  {\cal W}_j^{1,\infty;P_{+}}$, $g\in \widetilde{\cal W}_j^{1,\infty;P_{\pm}}([0,T])$.
 Assume either (i) $P_+(z)=e^{\lambda^{1/(\beta-1)}z}$, $P_-(z)=e^{c\lambda^{1/(\beta-1)}z}$ with $c>0$ smaller than some absolute constant ({\em exponential case}); 
 or (ii) $P_{\pm}(z)=z^{d_{\pm}}$
with $d_-\ge 1$, $d_+ - d_-> \frac{\beta-1}{\beta}d$. Then the   KPZ equation (\ref{eq:KPZdowntoj}) has a unique ${\cal W}_j^{1,\infty;P_-}$-solution.
}

\bigskip\noindent
{\bf Theorem 2 (estimates).}
{\em Under the conditions of Theorem 1, the following estimates hold for the solution $h$,
\BEQ ||| \lambda^{1/(\beta-1)} h_t|||_{{\cal W}_j^{1,\infty;P_-}}(x) \lesssim e^{-2^{-j}t}
P_-^{-1}\circ P_+\left( |||\lambda^{1/(\beta-1)} h_0|||_{{\cal W}^{1,\infty;P_+}}(x) \right) + 
P_-^{-1}\circ P_+\left(||| \lambda^{1/(\beta-1)} g|||_{\widetilde{\cal W}^{1,\infty;P_{\pm}}([0,t])}(x) \right)\EEQ
with: (i) $d'=d$ (exponential case); (ii) $d'=2(d_+-d_-)/d_-(\beta-1)$ (polynomial case). 
}

\bigskip\noindent
Comparing (\ref{intro:tilde-g-norm}) with (\ref{intro:def:g-norm}), one sees
there is an extra, rather innocuous multiplicative factor $(1+2^{-j}s)^{d'/2}$ which
may be absorbed by a suitable redefinition of the exponentially decreasing factor,
namely, $e^{-2^{-j}s} (1+2^{-j}s)^{d'/2}\le C e^{-c2^{-j}s}$ for a suitable constant
$C$ provided $c<1$.

\medskip
{\bf Remark.} 
 Regularized white noise $\eta$
belongs a.s. to all these functional spaces (as shown in the exponential case in  \cite{Unt-KPZ1}, section 6), and
one has an explicit  log-normal deviation formula for  its local time-dependent quasi-norms at a given space location $x$.


\subsection{Outline of the article}


The canvas is very simple. Generally speaking it is difficult to work directly
with (\ref{eq:KPZdowntoj}) with $u_0,g$ in ${\cal W}$-spaces since nothing guarantees a priori
the existence and unicity of the solution. So most of the work consists first in
establishing estimates in the case of "bounded data", i.e. for $u_0\in {\cal W}^{1,\infty}$, $g\in C([0,T],{\cal W}^{1,\infty})$. It is standard (from the general
theory of viscous Hamilton-Jacobi equations) that the equation has then a unique
solution in $C([0,T],{\cal W}^{1,\infty})$. It is also well-known that the equation
admits a Hamilton-Jacobi-Bellman reformulation, which is the object of section 1.
\S 2.1 is the main paragraph, dedicated to our estimates in local ${\cal W}$-quasi norms.
Then we show in \S 2.2 that the sequence $(h^{(n)})_{n\ge 1}$ of solutions of
the KPZ equations with bounded "cut-off" data $(h_0^{(n)},g^{(n)})$ converges locally
 uniformly to a unique function $h$, which can be shown to be a classical solution
 of the initial KPZ equation.
 
\medskip
\noindent Section 3 plays the r\^ole of an appendix. We collect some intermediate
results in \S 3.1. We sketch in \S 3.2 a derivation of estimates similar (but not
quite) to those
of Theorem 2 in terms of {\em poinwise} quasi-norms. The distinction between {\em poinwise} and {\em local} quasi-norms has been explained at length in our previous
article. Suffice it to say at this stage that, with a more cautious use of the $\locsup$ and
$\Locsup$ operations, functions in the associated {\em poinwise ${\cal W}-$spaces} are allowed to have arbitrarly
large local fluctuations, contrary to those in {\em local ${\cal W}$-spaces} (see
 Definition \ref{def:locsup} and ensuing discussion). The extension of our results to this much more general frame is satisfactory in itself, but it makes statements much more technical, which
is why we decided not to include this section in the main body of the article, nor
to write down all details.


\section{A Hamilton-Jacobi-Bellman reformulation of the  KPZ equation}


The general purpose of this section is to provide a random path representation, called
Hamilton-Jacobi-Bellman representation, for solutions of the KPZ or
infra-red cut-off KPZ equation. {\em We assume throughout that the data are bounded,}
i.e. $h_0\in {\cal W}^{1,\infty}, g\in C([0,T],{\cal W}^{1,\infty})$, in order to
apply the results of the theory in their standard form.


\subsection{The Hamilton-Jacobi-Bellman theory for viscous Hamilton-Jacobi equations}


For this paragraph which does not contain any new result, the reader may consult \cite{Fle}, \cite{Kry} or \cite{Sme}.
We restrict this introductory and somewhat loose discussion to a subclass of Hamilton-Jacobi equations, including
the "universality class" of the KPZ equation,
\BEQ \partial_t h=\Del h+ V(\nabla h)+g,\qquad h(t=0)=h_0 \label{eq:HJ} \EEQ
where $V:\R^d\to\R$ is assumed to be $C^2$ and convex, and $\frac{V(y)}{|y|}\to_{|y|\to \infty} +\infty$, so that $V$ has a well-defined
Legendre transform $\tilde{V}$ with the same properties; recall $\tilde{V}(\alpha):=\sup_{y\in\R^d}
(\alpha\cdot y-V(y)),\ \alpha\in\R^d$  (see e.g. \cite{Eva}, \S 3.3.2). 
Here we need more precisely $\widetilde{V(-\, \cdot)}(\cdot)=\tilde{V}(-\, \cdot)$.
The parameter conjugate to $y$ is $\nabla(V(-y))=-\nabla V(-y)$.  By
definition, $V$ and $\tilde{V}$ are related by  
\BEQ V(-y)=\sup_{\alpha} (\alpha y-\tilde{V}(-\alpha))=-\alpha^*(y)y-\tilde{V}(\alpha^*(y))\EEQ
where by definition  $\alpha^*(y)=+\nabla V(-y)$.

\medskip

Let  $B_t,\ t\ge 0$ be a $d$-dimensional Brownian motion, and consider the following class of stochastic differential
equations,
\BEQ dX^{\alpha}_s=\alpha_s ds+ dB_s, s\ge t \label{eq:diffusion} \EEQ
with initial condition $X^{\alpha}_t=x$, where $\alpha=(\alpha_s)_{s\ge 0}$ is an {\em admissible
strategy}, i.e. a progressively measurable, $\R^d$-valued
process with respect to the filtration defined by the Wiener process. We shall sometimes leave the dependence
of $X$ on $\alpha$ implicit and write $X$ instead of $X^{\alpha}$.
Fix a terminal time $T\ge 0$, a function $u_T:\R^d\to\R$  and some function $f=f(t,x;y)$ on $\R\times\R^d\times\R^d$. Bellman's original idea is to try and minimize the {\em cost functional}
 \BEQ J(t,x;\alpha):=\esper^{t,x}\left[ \int_t^T  f(s,X_s^{\alpha};\alpha_s)ds+u_T(X^{\alpha}_T)\right] \label{eq:cost} \EEQ
with respect to all admissible strategies. The notation $\esper^{t,x}[\cdot]$ emphasizes the initial condition
$X_t^{\alpha}=x$ for the diffusion. The result,
\BEQ v(t,x):=\inf_{\alpha} J(t,x;\alpha) \label{eq:inf-cost} \EEQ
is called the {\em value function}.

\medskip

Now, {\em Bellman's principle} states that, for $t\le \bar{t}\le T$,
\BEQ v(t,x)=\inf_{\alpha}\left( \esper^{t,x}\left[ \int_t^{\bar{t}} f(s,X_s^{\alpha};\alpha_s)ds + 
v(\bar{t},X^{\alpha}_{\bar{t}})\right] \right). \label{eq:Bellman's-principle}\EEQ
 This is essentially straightforward  since the choice of the optimal
strategy {\em after} time $\bar{t}$ depends by the Markov property of the Wiener process only on $X_{\bar{t}}$. Let now $\bar{t}=t+o(1)$ and
apply It\^o's formula. Note that the solution of (\ref{eq:Bellman's-principle}) is unique provided
one assumes the terminal condition $v_T=u_T$ (take $\bar{t}=T$). One gets, for $\bar{t}-t$
small, 
\BEQ \esper^{t,x}\left[ \int_t^{\bar{t}} f(s,X_s;\alpha_s)ds\right]=(\bar{t}-t)f(t,x;\alpha_t)+\ldots \EEQ
\BEQ \esper^{t,x}[v(\bar{t},X_{\bar{t}})]=v(t,x)+(\bar{t}-t)(\partial_t v(t,x)+{\cal L}^{\alpha}v(t,x))+\ldots \EEQ
where ${\cal L}^{\alpha}:=\Del+\alpha_t\cdot\nabla$ is the generator of the diffusion process (\ref{eq:diffusion}).
Taking the limit $\bar{t}\to t$ yields {\em{Bellman's differential equation}},
\BEQ \inf_{\alpha} \left\{ (\partial_t+{\cal L}^{\alpha})v(t,x)+f(t,x;\alpha_t) \right\}=0, \label{eq:Bellman} \EEQ
together with the obvious terminal condition, $v_T=u_T$.

\bigskip

Let us now choose
\BEQ f(t,x;y):=\tilde{V}(-y)-g(T-t,x)\EEQ
One immediately checks that Bellman's equation is equivalent to
\BEQ (\partial_t+\Del)v(t,x)+\inf_{\alpha}\left(\alpha_t\cdot\nabla v(t,x)+\tilde{V}(-\alpha_t)\right)-g(T-t,x)=0 \EEQ
or (letting $\alpha_t=-\alpha^*(\nabla v(t,x))$)
\BEQ (\partial_t+\Del)v(t,x)-V(\nabla v(t,x))-g(T-t,x)=0.\EEQ
Thus  one sees that 
\BEQ h(t,x):=-v(T-t,x), 0\le t\le T \label{eq:time-reversal} \EEQ
 satisfies the Hamilton-Jacobi equation (\ref{eq:HJ}).

\medskip

If the solution $h$ of the Hamilton-Jacobi equation (\ref{eq:HJ}) is unique, then the following
Feynman-Kac type formula holds.

\begin{Proposition} \label{prop:FK0}
Let 
\begin{itemize}
\item[(i)]
$h(t,\cdot)=-u(T-t,\cdot), \ t\ge 0$ solve the Hamilton-Jacobi equation (\ref{eq:HJ}) with initial
condition $h_0=-u_T$; 
\item[(ii)] $X^{\alpha^*}$ be the solution of the stochastic differential
\BEQ dX_s^{\alpha^*}=\alpha^*(\nabla u_s(X^{\alpha^*}_s))ds+dB_s,\qquad s\ge t\EEQ
with initial condition $X_t^{\alpha^*}=x$;
\item[(iii)] $v(t,x)$ be the function defined as
\BEQ v(t,x):=\esper^{t,x}\left[ \int_t^T \left( \tilde{V}(\alpha^*(\nabla u_s(X_s^{\alpha^*})))-g(T-s,X_s^{\alpha^*})\right) ds+u_T(X^{\alpha^*}_T)\right].
\label{eq:FK0}\EEQ
\end{itemize}
Then $v\equiv u$.
\end{Proposition}

{\bf Proof.} Clearly $v_T=u_T$.  We write $X=X^{\alpha^*}$ and let
\BEQ w_{t,x}(\bar{t}):=\esper^{t,x} \left[ \int_t^{\bar{t}} 
\left( \tilde{V}\left(\alpha^*(\nabla u_s(X_s))\right) - g(T-s,X_s) \right) ds +
u_{\bar{t}}(X_{\bar{t}})\right], \qquad t\le \bar{t}\le T.\EEQ  
  By It\^o's formula,
\BEA && \frac{d}{d\bar{t}} w_{t,x}(\bar{t})  \nonumber\\
&& = \esper^{t,x} \left[ \tilde{V}\left(\alpha^*(\nabla u_{\bar{t}}(X_{\bar{t}}))\right)  +\alpha^*(\nabla u_{\bar{t}}
(X_{\bar{t}}))\cdot \nabla u_{\bar{t}}(X_{\bar{t}}) + (\partial_{\bar{t}}+\Del)u_{\bar{t}}(X_{\bar{t}}) -g(T-\bar{t},X_{\bar{t}}) \right].
\label{eq:FK0-Ito} \EEA
The first two summands in (\ref{eq:FK0-Ito}) sum up by definition of $\alpha^*$ to $-V(-\nabla u_{\bar{t}}(X_{\bar{t}}))$.
Since $-u(T-t,\cdot)$ solves eq. (\ref{eq:HJ}), $\frac{d}{d\bar{t}} w_{t,x}(\bar{t})=0$.  Thus
$w_{t,x}(T)=w_{t,x}(t)$, or in other terms, $v=u$. \hfill\eop

\medskip

Summarizing the above discussion:

\begin{Proposition} \label{prop:FK1}
Assume the solution of the Hamilton-Jacobi equation (\ref{eq:HJ}) is unique. Then
\BEA  h(t,x) &=& \sup_{\alpha} \esper^{0,x}\left[ \int_0^t \left(-\tilde{V}(-\alpha_s)+g(t-s,X_s^{\alpha})
\right) ds + h_0(X_t^{\alpha}) \right] \nonumber\\
&=& \esper^{0,x}\left[ \int_0^t \left(-\tilde{V}(-\alpha^*_s)+g(t-s,X_s^{\alpha^*})
\right) ds + h_0(X_t^{\alpha^*}) \right] 
\EEA
where $\alpha_s^*:=\alpha^*(-\nabla h_{t-s}(X_s^{\alpha^*}))$.
\end{Proposition}

\medskip


\subsection{Application to the infra-red cut-off KPZ equation}

We now consider the infra-red cut-off KPZ equation (\ref{eq:KPZdowntoj}), with the supplementary assumption
$\frac{V(y)}{y}\to_{y\to\infty}+\infty$ (this supplementary assumption is necessary to define the Legendre
transform of $V$, but we show in section 2 how to get rid of it very easily). Compared to (\ref{eq:HJ}), $V$ has been
changed to $\lambda V$, and an extra linear term $-2^{-j}h$ appears in the right-hand side. This accounts for two minor modifications with respect to the above analysis. First,
a simple scaling argument yields 
$\widetilde{\lambda V}(p)=\lambda \tilde{V}(\frac{p}{\lambda})$.  Second, the modified generator
${\cal L}^{\alpha}_{killed}(t,x):= \Del +\alpha_t \cdot\nabla -2^{-j}$ corresponds to a diffusion with 
killing rate $2^{-j}$.  A straightforward extension of the results of the previous paragraph yields

\begin{Lemma}[Hamilton-Jacobi-Bellman representation of the infra-red cut-off KPZ equation] \label{lem:HJB-representation}
Let $h$ be the solution of (\ref{eq:KPZdowntoj}). Then 
\BEQ h(t,x)=\sup_{\alpha} J_{\alpha}(t,x), \label{eq:HJB-representation}\EEQ where
\BEQ J_{\alpha}(t,x):=\esper^{0,x} \left[  \int_0^t e^{-2^{-j}s} \left(-\lambda \tilde{V}(\frac{\alpha_s}{\lambda}) 
+g(t-s,X^{\alpha}_s) \right) ds + e^{-2^{-j} t} h_0(X^{\alpha}_t)\right]  \label{eq:Jalpha} \EEQ
where $X^{\alpha}$ is the solution of the stochastic differential equation, 
\BEQ dX^{\alpha}_s=\alpha_s ds+ dB_s \label{eq:SDE} \EEQ
 with
initial condition $X^{\alpha}_0=x$.  An explicit optimal path, $X^{\alpha^*}$,  is given by $\alpha^*_s=\lambda \nabla V(-\nabla h_{t-s}(X_s^{\alpha^*}))$.
\end{Lemma}

{\bf Remark.} 
Let us give an elementary application of Lemma \ref{lem:HJB-representation}. Taking $\alpha\equiv 0$ in eq. (\ref{eq:Jalpha}), one gets
\BEQ h(t,x)\ge e^{-2^{-j}t}\inf h_0 - \int_0^t e^{-2^{-j}s} ||g_{t-s}||_{\infty} ds\EEQ
since $\tilde{V}(0)=\max_p (-V(p))=0$.
On the other hand, $\tilde{V}(p)\ge -V(0)=0$ so
 \BEQ h(t,x)\le e^{-2^{-j}t}\sup h_0 + \int_0^t e^{-2^{-j}s} ||g_{t-s}||_{\infty} ds.\EEQ
 Of course, both inequalities also follow from a direct application of the maximum principle.
 
 \bigskip


\section{PDE estimates}


We prove in this section our main results, Theorem 1 and Theorem 2 of the Introduction. The strategy is roughly as follows. We first show that the estimates in local ${\cal W}$-quasi norms of Theorem 2 hold for {\em bounded} data, i.e. for $u_0\in {\cal W}^{1,\infty}, g\in C([0,T],{\cal W}^{1,\infty})$ (\S 2.1). Then we show that all sequences $h^{(n)}$ constructed as in Definition \ref{intro:def:Wj1inftyPsolution} converge
to a unique limit, which is a classical solution of (\ref{eq:KPZdowntoj}) satisfying
also the estimates of Theorem 2.


\subsection{A priori estimates}


{\em In this paragraph, we assume bounded data}, i.e. $u_0\in {\cal W}^{1,\infty}, g\in C([0,T], {\cal W}^{1,\infty})$, and show the estimates of Theorem 2. Time-dependent 
${\cal W}$-spaces
$\widetilde{\cal W}^{1,\infty;P_{\pm}}_j(x)$ depend implicitly on an extra index which we call
throughout 
$d'$. 

\medskip

We first bound  local ${\cal W}$-quasi norms of the solution itself, then we take
care of those of its gradient. To get an a priori bound of  the solution alone, we do not need to have a control over the gradient of the data $h_0,g$. Thus it is more natural
to express bounds in terms of the quasi-norms $|||\, \cdot\,  |||_{{\cal W}^{0,\infty;P}}(x)$ and $|||\, \cdot\,  |||_{\widetilde{\cal W}^{0,\infty;P_{\pm}}([0,t])}(x)$ obtained from (\ref{eq:Wj1inftyP}), (\ref{intro:tilde-g-norm})
by leaving out  the gradient terms. In any case, clearly
$|||\, \cdot\,  |||_{{\cal W}^{0,\infty;P}}(x)\le |||\,  \cdot\, |||_{{\cal W}^{1,\infty;P}}(x)$
and $|||\,  \cdot\, |||_{\widetilde{\cal W}^{0,\infty;P_{\pm}}([0,t])}(x) \le |||\, \cdot\, |||_{\widetilde{\cal W}^{1,\infty;P_{\pm}}([0,t])}(x)$.

\medskip \noindent Before we state our main result, we make a small remark.
The function $h$ is obtained as the solution of some {\em maximization} problem,
$h(t,x)=\sup_{\alpha} J_{\alpha}(t,x)$. Thus $h(t,x)\ge \phi(t,x):=J_0(t,x)$, where
$\phi$ is the solution of the linearized problem, $(\partial_t-\Del+2^{-j})\phi(t,x)=g$,
with initial condition $\phi_0\equiv h_0$.
Hence (by the comparison principle), $-h\le \psi$, where $\psi$ is the solution of the KPZ equation
$(\partial_t-\Del+2^{-j})\psi=\lambda V(\nabla\psi)-g$ with initial condition
$\psi_0\equiv -h_0$. Thus we have a two-sided bound for $h$, $|h(t,x)|\le \sup_{\alpha} \tilde{J}_{\alpha}(t,x)$, with
\BEQ  \tilde{J}_{\alpha}(t,x):=\esper^{0,x} \left[  \int_0^t e^{-2^{-j}s} \left(-\lambda \tilde{V}(\frac{\alpha_s}{\lambda}) 
+|g|(t-s,X^{\alpha}_s) \right) ds + e^{-2^{-j} t} |h_0|(X^{\alpha}_t)\right].  \label{eq:Jalphatilde} \EEQ

\begin{Theorem} \label{th:main}
\begin{itemize}
\item[(i)] (exponential case)
It holds for $d'=d$
\BEQ \ln (e^{\lambda^{1/(\beta-1)} |h_t|})^*(x) \lesssim \lambda^{1/(\beta-1)} \left( 
e^{-2^{-j} t} |||\lambda^{1/(\beta-1)} h_0|||_{{\cal W}_j^{0,\infty;\lambda^{1/(\beta-1)}}}(x) + |||\lambda^{1/(\beta-1)} g|||_{\widetilde{\cal W}_j^{0,\infty;
\lambda^{1/(\beta-1)}}([0,T])}(x)  \right).  \label{eq:main} \EEQ

\item[(ii)] Let  $d_-\ge 1$, $d_+>d_- +\frac{\beta-1}{\beta} d$, and $P_{\pm}(z)=z^{d_{\pm}}$. Then it holds for $d'=2(d_+-d_-)/d_-(\beta-1)$
\BEA &&  P_-^{-1}\left(\left( P_-(\lambda^{1/(\beta-1)} |h_t|)\right)^*(x)\right)\lesssim e^{-2^{-j}t}  P_-^{-1}\circ P_+ \left( |||\lambda^{1/(\beta-1)} h_0|||_{{\cal W}_j^{0,\infty;
P_+}}(x)\right) \nonumber\\
&& \qquad \qquad \qquad  + P_-^{-1}\circ P_+\left( |||\lambda^{1/(\beta-1)}g|||_{\widetilde{\cal W}_j^{0,\infty;P_{\pm}}([0,T])}(x)\right). \label{eq:main-pol} \EEA

\end{itemize}
\end{Theorem}

\noindent Note that if $\beta=2$, (\ref{eq:main}) is a less precise statement of 
the bound found in \cite{Unt-KPZ1}, Theorem 2. Note that (\ref{eq:main}) has
been written in this way for simplicity and to make comparison to \cite{Unt-KPZ1}
easier, but actually it may also equivalently be stated as (\ref{eq:main-pol}), with
$P_-(z)=P_+(z)=e^{\lambda^{1/(\beta-1)}z}$. On the other hand, in the polynomial case (ii), (\ref{eq:main-pol}) reads more explicitly
\BEA &&  \left( (|h_t|^{d_-})^*(x) \right)^{1/d_-} \lesssim \lambda^{(\frac{d_+}{d_-}-1)/(\beta-1)}  \left(  e^{-2^{-j}t} 
 \left\{ \left( |h_0|^{d_+})^*(x) \right)^{1/d_+} \right\}^{d_+/d_-} + \right.\nonumber\\
 && \qquad\qquad \qquad \left. +
 \int_0^t e^{-2^{-j}s} (1+2^{-j}s)^{d'/2} \left\{ \left( (2^j |g(t-s,\cdot)|)^{d_+})^*(x) \right)^{1/d_+} \right\}^{d_+/d_-} \right). \EEA

\bigskip

\noindent {\bf Proof.} 
We give the proof in the case when $V(y)/y\to_{y\to\infty} +\infty$. Otherwise one should simply note that a uniform
bound is obtained for $V$ replaced by $V_{\eps}(y):=V(y)+\eps y^2$ and let $\eps\to 0$. 
 
\medskip

\noindent (i) ({\em exponential case}). 
Recall that $|h(t,x)|\le \sup_{\alpha} \tilde{J}_{\alpha}(t,x)$, where $\tilde{J}_{\alpha}$ is defined in (\ref{eq:Jalphatilde}).
By Jensen's inequality, 
\BEQ e^{\lambda^{1/(\beta-1)} |h|(t,x)}\le \sup_{\alpha} J_{\alpha}(\lambda;t,x),\EEQ where
\BEQ J_{\alpha}(\lambda;t,x):=\esper^{0,x} \left[ \exp \lambda^{1/(\beta-1)} \left( \int_0^t e^{-2^{-j}s} \left(-\lambda \tilde{V}(\frac{\alpha_s}{\lambda}) 
+|g|(t-s,X_s) \right) ds + e^{-2^{-j} t} |h_0|(X_t) \right) \right] \label{eq:Jlambda} \EEQ 
with $X_t=X_t^{\alpha}$.
The proof will be divided into three steps.

\begin{itemize}
\item[1.] 
Since
\BEA  \int_0^t \tilde{V}(\frac{\alpha_s}{\lambda}) e^{-2^{-j} s} ds & \ge& 2^{-j} \int_0^t \tilde{V}(\frac{\alpha_s}{\lambda}) \ \cdot\ \left( \int_s^t e^{-2^{-j}u} du\right) ds \nonumber\\
&=& 2^{-j} \int_0^t \left( \int_0^u \tilde{V}(\frac{\alpha_s}{\lambda})ds\right) e^{-2^{-j}u} du,\EEA
one gets
\BEA &&  J_{\alpha}(\lambda;t,x)\le \limsup_{n\to\infty} \esper^{0,x} \exp \left( \sum_{k=1}^n e^{-2^{-j} kt/n}
\left\{ -\half 2^{-j}\frac{t}{n} \lambda^{\beta/(\beta-1)} \int_0^{kt/n} \tilde{V}(\frac{\alpha_s}{\lambda})
ds+  \frac{t}{n} \lambda^{1/(\beta-1)} |g_{t-kt/n}|(X_{kt/n}) \right\}  \right. \nonumber\\
&& \qquad \qquad \qquad \qquad \left. +e^{-2^{-j} t} \left\{ -\half \lambda^{\beta/(\beta-1)} \int_0^t
\tilde{V}(\frac{\alpha_s}{\lambda})ds+\lambda^{1/(\beta-1)} |h_0|(X_t) \right\} \right).\EEA
Let $p:=(2^{-j}\frac{t}{n})^{-1}-1$ (see \cite{Unt-KPZ1}, proof of Lemma 4.4) and
$p_k:=\frac{1}{p+1}\left(\frac{p}{p+1}\right)^k$ ($k=0,\ldots,n-1$), $p_n:=\left(\frac{p}{p+1}\right)^n$, so that $\frac{1}{p_0}+\ldots+\frac{1}{p_n}=1$ and, for $n$ large, 
\BEQ p_k\sim 2^{-j}\frac{t}{n} e^{-2^{-j} kt/n}, \qquad p_n\sim e^{-2^{-j}t}.\EEQ
 Using the generalized H\"older property,
$$\esper[Y_0\ldots Y_n]\le \prod_{k=0}^n \esper[Y_k^{p_k}]^{1/p_k},  \qquad (Y_0,\ldots,Y_n\ge 0),$$
one gets
\BEA && J_{\alpha}(\lambda;t,x) \le  \left( \esper^{0,x}\left[
 \exp\left(-\half \lambda^{\beta/(\beta-1)} \int_0^t \tilde{V}(\frac{\alpha_s}{\lambda})ds
 +\lambda^{1/(\beta-1)} |h_0|(X_t) \right)\right] \right)^{ e^{-2^{-j}t}}  \nonumber\\
&& \qquad \limsup_{n\to\infty} \prod_{k=1}^n \left( \esper^{0,x}\left[
\exp \left(-\half\lambda^{\beta/(\beta-1)}  \int_0^{kt/n} \tilde{V}(\frac{\alpha_s}{\lambda})
ds  + \lambda^{1/(\beta-1)} 2^j |g_{t-kt/n}|(X_{kt/n}) \right) \right] \right)^{ 2^{-j}\frac{t}{n}
e^{-2^{-j}kt/n}}. \nonumber\\ \label{eq:3.5}
\EEA

\item[2.] We shall bound each individual term 
\BEQ \bar{J}(\lambda;u,x):=\esper^{0,x}\left[
\exp \left(-\half \lambda^{\beta/(\beta-1)}  \int_0^{u} \tilde{V}(\frac{\alpha_s}{\lambda})
ds+ \lambda^{1/(\beta-1)} 2^j|g_{t-u}|(X_{u}) \right) \right], \label{eq:3.6} \EEQ
$u=kt/n$ $(k=1,\ldots,n)$  
in the above expression. The factor depending on the initial condition is identical to $\bar{J}(\lambda;t,x)$
except that $2^j|g_t|$ is replaced with $|h_0|$, so we do not discuss it any more and assume $h_0\equiv 0$ in the
sequel to simplify notations. Note also that the generalized H\"older property used in the previous paragraph may also be applied to $e^{\tau\Del}\esper\left[\ \cdot\ 
\right]$; thus $(e^{\tau\Del}J_{\alpha}(\lambda;t))(x)$ is bounded by the weighted product of the $(e^{\tau\Del}\bar{J}(\lambda;u))(x)$.

Note first that
\BEQ \tilde{V}(p)=\sup_x (p\cdot x-V(x))\ge (|p|^2-V(p)) {\bf 1}_{|p|\le 1}+ \left(|p|^{\beta/(\beta-1)}-V(|p|^{1/(\beta-1)}) \right) {\bf 1}_{|p|>1}\ge \half \min(|p|^2,|p|^{\beta/(\beta-1)}) \EEQ
hence 
\BEQ \widetilde{\lambda V}(p)=\lambda \tilde{V}(\frac{p}{\lambda})\ge\half\min\left(\frac{|p|^2}{\lambda},
\left( \frac{|p|^{\beta}}{\lambda}\right)^{1/(\beta-1)} \right).\EEQ

Let $\Omega:=\{s\in[0,u];\ |\alpha_s|\le \lambda\}$, and $\bar{\Omega}:=[0,u]\setminus\Omega$. On $\Omega$,
$\lambda^2 \tilde{V}(\frac{\alpha_s}{\lambda})\ge \frac{\alpha_s^2}{2}$; on $\bar{\Omega}$, 
$\lambda^{\beta/(\beta-1)} \tilde{V}(\frac{\alpha_s}{\lambda})\ge \frac{|\alpha_s|^{\beta/(\beta-1)}}{2}$. We now
distinguish two cases. If $\int_{\Omega}|\alpha_s|ds\ge \half \int_0^u |\alpha_s|ds$, then
\BEQ \lambda^{\beta/(\beta-1)} \int_0^u \tilde{V}(\frac{\alpha_s}{\lambda})ds\ge \lambda^2 \int_{\Omega}
\tilde{V}(\frac{\alpha_s}{\lambda})ds\ge \half \int_{\Omega}\alpha_s^2 ds\ge \frac{1}{2|\Omega|}
(\int_{\Omega}|\alpha_s|ds)^2\ge \frac{1}{8u} (\int_0^u \alpha_s ds)^2.\EEQ 
 Otherwise it follows from H\"older's
inequality that
\BEA && \lambda^{\beta/(\beta-1)} \int_0^u \tilde{V}(\frac{\alpha_s}{\lambda})ds\ge \half \int_{\bar{\Omega}}
|\alpha_s|^{\beta/(\beta-1)} ds\ge \half |\bar{\Omega}|^{-1/(\beta-1)} (\int_{\bar{\Omega}} |\alpha_s|ds)^{\beta/(\beta-1)} \nonumber\\
&&\qquad \qquad \qquad \qquad \ge \half 2^{-\beta/(\beta-1)} u^{-1/(\beta-1)} \left|\int_0^u \alpha_s ds \right|^{\beta/(\beta-1)}.\EEA

 By definition, $X_u-x=\int_0^u \alpha_s ds +  B_u$, hence the net outcome of all these computations is
\BEQ \half\lambda^{\beta/(\beta-1)} \int_0^u \tilde{V}(\frac{\alpha_s}{\lambda})ds\ge  \frac{1}{64} 
\min\left( \frac{|X_u-x|^2}{u}, \left(\frac{|X_u-x|^{\beta}}{u}\right)^{1/(\beta-1)} \right) {\bf 1}_{|X_u-x|\ge
2|B_u|}.\EEQ

  Hence $\bar{J}(\lambda;u,x) \le  1+ \sup_{B(x,2^{j/2})} (e^{\lambda^{1/(\beta-1)} 
  2^j|g_{t-u}|}-1) + I_1(x)+I_2(x)+I_3(x)$, where:
\BEA && I_1(x):=   \esper^{0,x} \left[ {\bf 1}_{2^{j/2}< |X_u-x|<2|B_u|} 
(e^{\lambda^{1/(\beta-1)} 2^j|g_{t-u}|(X_u)} -1) \right], \nonumber\\
&& \qquad \qquad 
I_2(x):= \esper^{0,x} \left[ {\bf 1}_{|X_u-x|\ge \max(2^{j/2}, 2|B_u|)} e^{-\frac{|X_u-x|^2}{64u}} (e^{\lambda^{1/(\beta-1)}  2^j|g_{t-u}|(X_u)}-1)  \right], \qquad \nonumber\\
&& \qquad\qquad\qquad
I_3(x):=
 \esper^{0,x} \left[ {\bf 1}_{|X_u-x| \ge \max(2^{j/2},2|B_u|)} e^{-\frac{1}{64} (|X_u-x|^{\beta}/u)^{1/(\beta-1)}} (e^{\lambda^{1/(\beta-1)}2^j|g_{t-u}|(X_u)}-1) \right] 
 \nonumber\\ 
 \label{eq:3.13} \EEA

\medskip

Estimates of $I_1,I_2,I_3$ are all based on Lemma \ref{lem:local-integration}.
We  first deal with $I_1$. Let $\tilde{\Del}\in \widetilde{\D}^j$ be a cube of size $c2^{j/2}$ (see proof of Lemma \ref{lem:local-integration} for notations). then 
$$\proba^{0,x}\left[ \left(|X_u-x|<
2|B_u| \right)\cap\left(X_u-x\in\tilde{\Del} \right)\right]\le\proba^{0,x}
\left[ |B_u|>\half \min_{y\in\tilde{\Del}} |y| \right] \lesssim e^{-c'(\min_{y\in\tilde{\Del}} |y|)^2/u}$$ for
some constant $c'>0$.   Thus the Lemma asserts that, for every $\tau\ge 0$,
\BEQ  e^{\tau\Del}I_1(x)  \lesssim  (2^{-j}u)^{d/2} \ 
\locsup^j\left( e^{\lambda^{1/(\beta-1)}2^j|g_{t-u}| }-1 \right)^*(x).\EEQ

For $I_2$ and $I_3$ we  bound $\proba[X_u-x\in\Del]$ by $1$ and use instead the
 exponential decreasing weights $e^{-\frac{|X_u-x|^2}{64u}}\approx u^{d/2}
 \Phi_{64u}^2(X_u-x)$, resp. $e^{-\frac{1}{64} (|X_u-x|^{\beta}/u)^{1/(\beta-1)}}\approx
u^{d/\beta} \Phi_{64u}^{\beta}(X_u-x)$ (see \S 3.1 for notations), yielding
\BEQ  e^{\tau\Del}I_2(x)  \lesssim  (2^{-j}u)^{d/2} \ 
\locsup^j\left( e^{\lambda^{1/(\beta-1)} 2^j|g_{t-u}| }-1 \right)^*(x)\EEQ
and
\BEQ  e^{\tau\Del} I_3(x)  \lesssim  (2^{-j})^{d/2} u^{d/\beta} \ 
\locsup^j\left( e^{\lambda^{1/(\beta-1)} 2^j|g_{t-u}| }-1 \right)^*(x).\EEQ

\medskip
To finish with, we note that
\BEQ \locsup^j\left( e^{\lambda^{1/(\beta-1)} 2^j|g_{t-u}|}-1 \right)^*(x) =
\locsup^j\left( e^{\lambda^{1/(\beta-1)}2^j|g_{t-u}|} \right)^*(x)-1=e^{\lambda^{1/(\beta-1)} ||2^j|g_{t-u}|||_{{\cal W}^{0,\infty;\lambda^{1/(\beta-1)}}} (x)} - 1 \EEQ
hence
\BEQ (e^{\tau\Del} \bar{J}(\lambda;u))(x)\le 1+ C(1+(2^{-j}u)^{d/2}) \left\{ \exp\left(\lambda^{1/(\beta-1)} ||| 2^j g_{t-u}|||_{
{\cal W}_j^{0,\infty;\lambda^{1/(\beta-1)}}} (x)\right) -1 \right\}  \label{eq:3.27} \EEQ 
for some constant $C>0$. 

\item[3.]
Combining (\ref{eq:3.5}), (\ref{eq:3.6}) and (\ref{eq:3.27}), one gets
\BEQ \ln (J_{\alpha}(\lambda;t))^*(x)\lesssim \limsup_{n\to\infty} 2^{-j} \frac{t}{n} \sum_{k=1}^n e^{-2^{-j} kt/n}
\ln\left(1+C(1+(2^{-j}kt/n)^{d/2})f(kt/n) \right), \label{eq:lnJalpha} \EEQ
where $f(u):=\exp\left(\lambda^{1/(\beta-1)} ||2^j g_{t-u} |||_{{\cal W}_j^{0,\infty;\lambda^{1/(\beta-1)}}}(x)\right) -1$.

By definition, $ 2^{-j} \frac{t}{n} \sum_{k=1}^n e^{-2^{-j} kt/n} \ln(1+f(kt/n))\to_{n\to\infty} \lambda^{1/(\beta-1)}
|||g|||_{{\cal W}_j^{0,\infty;\lambda^{1/(\beta-1)}}([0,t])}(x).$ We must still take into account the polynomial correction
$C(1+(2^{-j}u)^{d/2})$ in factor of $f(u)$ in the right-hand side of (\ref{eq:lnJalpha}). Let  
$\Omega_1:=\{u\in[0,t]\ |\ f(u)\gg 1\}$ and $\Omega_2:=[0,t]\setminus\Omega_1$. On $\Omega_1$,
$e^{-2^{-j}u}\ln(1+C(1+(2^{-j}u)^{d/2})f(u))\lesssim e^{-2^{-j}u}(\ln(1+(2^{-j}u)^{d/2})+\ln(1+f(u)))\lesssim (1+2^{-j}u)^{d/2}e^{-2^{-j}u}\ln(1+f(u)).$
Similarly, on $\Omega_2$, $e^{-2^{-j}u}\ln(1+C(1+(2^{-j}u)^{d/2})f(u))\lesssim (1+2^{-j}u)^{d/2} e^{-2^{-j}u} f(u)\lesssim (1+2^{-j}u)^{d/2} e^{-2^{-j}u}\ln(1+f(u)).$ 
Thus the polynomial correction may be replaced by  the multiplicative factor $(1+2^{-j}u)^{d/2}$, which is precisely the factor 
included into the definition of  $|||\ \cdot \ |||_{\widetilde{\cal W}_j^{0,\infty;\lambda^{1/(\beta-1)}}([0,t])}(x)$.

\end{itemize}

\medskip

\noindent (ii) ({\em polynomial case}). 
The proof is very similar and we only emphasize the differences. We start from $\tilde{J}_{\alpha}(t,x)$ as in (i) and apply
the following inequality,  consequence of the generalized H\"older property for real-valued random variables $(Y_k)_k$,
\BEQ \esper\left[ \left(\sum_k \frac{1}{p_k} (Y_k)_+ \right)^{d_-}\right]  \le \left( \sum_k \frac{1}{p_k}\, \left(\esper[(Y_k)_+^{d_-}]\right)^{1/d_-} \right)^{d_-} \label{eq:muk} \EEQ
(expand in both sides, and bound each term $\esper[(Y_{k_1})_+\cdots (Y_{k_{d_-})_+}]$
in the left-hand side by $\prod_{l=1}^{d_-} \esper[(Y_{k_l})_+^{d_-}]^{1/d_-}$),
with
 $Y_k=\lambda^{1/(\beta-1)}
(-\half \lambda \int_0^{kt/n} \tilde{V}(\frac{\alpha_s}{\lambda}) ds +2^j |g_{t-kt/n}|(X_{kt/n}))$. Hence 
 (compare with  (\ref{eq:3.5}) and (\ref{eq:3.6})), extending $P_-$ to zero on $\R_-$,
\BEA P_-(\esper[\lambda^{1/(\beta-1)} \tilde{J}_{\alpha}(t,x)]) \le
\esper[ P_-(\lambda^{1/(\beta-1)} \tilde{J}_{\alpha}(t,x))] \lesssim  P_-\left(  \limsup_{n\to\infty} 2^{-j}\frac{t}{n}
\sum_{k=1}^n e^{-2^{-j}kt/n} P_-^{-1}(\bar{J}(P_-;kt/n,x)) \right). \nonumber\\
 \EEA
Instead of (\ref{eq:3.13}), one obtains
\BEA &&  \bar{J}(P_-;u,x)\le \esper^{0,x}\left[ {\bf 1}_{|X_u-x|<2|B_u|} P_-(\lambda^{1/(\beta-1)} 2^j|g_{t-u}|(X_u)) \right]
+ \esper^{0,x} \left[ {\bf 1}_{|X_u-x|\ge 2|B_u|} P_-\left(\lambda^{1/(\beta-1)}2^j|g_{t-u}|(X_u)-
\frac{|X_u-x|^2}{64u} \right)
\right]  \nonumber\\
&& \qquad \qquad \qquad + \esper^{0,x} \left[ {\bf 1}_{|X_u-x|\ge 2|B_u|} P_-\left(\lambda^{1/(\beta-1)} 2^j|g_{t-u}|(X_u)-
\frac{(|X_u-x|^{\beta}/u)^{1/(\beta-1)}}{64} \right)
\right] \nonumber\\
&& \qquad \qquad \qquad  \le \sup_{B(x,2^{j/2})} P_-(\lambda^{1/(\beta-1)} 2^j|g_{t-u}|(X_u))+ I_1(x)+I_2(x)+I_3(x)  \nonumber\\
 \EEA
where $I_1,I_2,I_3$ are analogous to the homonymous quantities considered in  the exponential case (i),  and consequently, $e^{\tau\Del}I_1(x)\lesssim (2^{-j}u)^{d/2}
\locsup^j \left(P_-(\lambda^{1/(\beta-1)} 2^j |g_{t-u}|) \right)^*(x).$ 

\medskip\noindent
 After these rather general
considerations, we now fix $P_{\pm}(z)=z^{d_{\pm}}$ with $d_+\ge d_-\ge 1$, and estimate $I_2$, $I_3$ by some specific arguments. We have $  P_-\left( \lambda^{1/(\beta-1)} 2^j|g_{t-u}|(x+y) -\frac{|y|^2}{64 u} \right)=0$ if $\lambda^{1/(\beta-1)}2^j|g_{t-u}|(x+y)\le |y|^2/64 u$, otherwise (using the trivial bound $(a-b)^{d_-}\le \frac{a^{d_+}}{b^{d_+-d_-}}$, valid for 
$a\ge b>0$), 
\BEQ  P_-\left( \lambda^{1/(\beta-1)} 2^j|g_{t-u}|(x+y) -\frac{|y|^2}{64 u} \right)
\le  \frac{\left(\lambda^{1/(\beta-1)} 2^j|g_{t-u}|(x+y)\right)^{d_+}}{(|y|^2/64u)^{d_+-d_-}}.\EEQ
 The  kernel
$\Phi(y)=(2^{j/2})^{2(d_+-d_-)-d} (|y|^2)^{-(d_+-d_-)}{\bf 1}_{|y|>2^{j/2}}$ is integrable at infinity since (by hypothesis)
$d_+-d_->\frac{\beta-1}{\beta}d\ge d/2$, and satisfies the hypotheses of Lemma \ref{lem:beta} (up to a normalization constant of order 1). Hence, by Lemma \ref{lem:local-integration},
\BEQ  e^{\tau\Del}I_2 (x) \lesssim  (2^{-j}u)^{d_+-d_-} 
P_+\left( |||\lambda^{1/(\beta-1)} 2^j|g_{t-u}|\  |||_{{\cal W}_j^{0,\infty;P_+}}(x)\right). \EEQ

\medskip
The same arguments apply to $I_3$, except that the replacement of $|y|^2/u$ by
$(|y|^{\beta}/u)^{1/(\beta-1)}$ leads to a normalized kernel, $\bar{\Phi}(y)\approx
(2^{j/2})^{\frac{\beta}{\beta-1}(d_+-d_-)-d} |y|^{-\frac{\beta}{\beta-1} (d_+-d_-)} {\bf 1}_{|y|>2^{j/2}}$  (we use here the integrability hypothesis $d_+-d_->\frac{\beta-1}{\beta} d$).
Thus
\BEA  e^{\tau\Del}I_3(x)  &\lesssim & (2^{j/2})^{-(d_+-d_-)\beta/(\beta-1)} u^{(d_+-d_-)/(\beta-1)} P_+\left( |||\lambda^{1/(\beta-1)} 2^j|g_{t-u}|\ |||_{{\cal W}_j^{0,\infty;P_+}}(x)\right) \nonumber\\
&\lesssim&
 (2^{-j}u)^{(d_+-d_-)/(\beta-1)} P_+\left( |||\lambda^{1/(\beta-1)} 2^j|g_{t-u}| \ |||_{{\cal W}_j^{0,\infty;P_+}}(x)\right). \EEA
The polynomial correction is in this case $P_-^{-1}\left((1+2^{-j}u)^{(d_+-d_-)/(\beta-1)} \right)$, whence $d'=2(d_+-d_-)/d_-(\beta-1)$.
\hfill \eop

\bigskip

An analogous bound on the gradient of
the solution is easily derived from our Theorem.

\begin{Corollary}[gradient bound] \label{cor:main}
Let $h_0\in {\cal W}_j^{1,\infty;P_+}$,  $g\in {\cal W}_j^{1,\infty;P_{\pm}}([0,T])$ and 
$(P_-,P_+)$ a pair of functions satisfying the hypotheses of Theorem 1. Then $h_t\in {\cal W}_j^{1,\infty;P_-}$ for all $t\le T$, and 
\BEQ  P_-^{-1}\left(\left( P_-(\lambda^{1/(\beta-1)} |2^{j/2}\nabla h_t|)\right)^*(x)\right) \lesssim 
e^{-2^{-j} t}  P_-^{-1}\circ P_+\left( |||\lambda^{1/(\beta-1)}h_0|||_{{\cal W}_j^{1,\infty;P_+}}(x)\right) +  P_-^{-1}\circ P_+\left( ||| \lambda^{1/(\beta-1)}g|||_{\widetilde{\cal W}_j^{1,\infty;
P_{\pm}}([0,T])}(x)\right). \label{eq:3.30} \EEQ
\end{Corollary}

{\bf Proof.} 

\begin{itemize}
\item[(i)] (gradient of the solution)
Let $\alpha^*$ maximize the right-hand side of (\ref{eq:HJB-representation}), $\vec{\eps}\in B(0,1)\setminus\{0\}$, and 
$\tilde{\del}_{\vec{\eps}}h(t,x):=\frac{h(t,x+\vec{\eps})-(1-2^{-j/2}|\vec{\eps}|)h(t,x)}{|\vec{\eps}|}$, $\tilde{\del}_{\vec{\eps}}g(t,x):=\frac{g(t,x+\vec{\eps})-(1-2^{-j/2}|\vec{\eps}|)g(t,x)}{|\vec{\eps}|}$ (see \cite{Unt-KPZ1}, section 4.3 for a similar proof). Then
\BEQ \tilde{\del}_{\vec{\eps}}h(t,x)\le \esper^{0,x} \left[ \int_0^t e^{-2^{-j}s} \left( -2^{-j/2} \lambda \tilde{V}(\frac{\alpha^*_s}{\lambda})
+ \tilde{\del}_{\vec{\eps}}g(t-s,X_s^{\alpha^*}) \right) ds+e^{-2^{-j}t} \tilde\del_{\vec{\eps}}h_0(X^{\alpha^*}_t). \right] \EEQ
Proceeding as in the proof of Theorem \ref{th:main} and using 
$$\tilde{\del}_{\vec{\eps}}h(t,\cdot)=\frac{h(t,\cdot+\vec{\eps})-
h(t,\cdot)}{|\vec{\eps}|}+2^{-j/2} h(t,\cdot)=-(1-|\vec{\eps}|) \tilde{\del}_{-\vec{\eps}}^j h(t,\cdot + 2^{j/2} \vec{\eps})+ (2-|\vec{\eps}|)h(t,\cdot + 2^{j/2} \vec{\eps})$$
(implying a two-sided bound, see \cite{Unt-KPZ1}, proof of Lemma 3.12),
 $\vec{\del}_{\vec{\eps}}g(t,\cdot)=\frac{g(t,\cdot+\eps)-
g(t,\cdot)}{|\eps|}-2^{-j/2} g(t,\cdot)$, one gets (\ref{eq:3.30}). In the
exponential case, H\"older's inequality, $e^{\tau\Del}(e^{\lambda(f_1+f_2)})\le 
\left( e^{\tau\Del}(e^{pf_1}) \right)^{1/p} \left( e^{\tau\Del} (e^{qf_2}) \right)^{1/q}$,
$p,q\ge 1$, $\frac{1}{p}+\frac{1}{q}=1$ must be used several times, implying the loss of regularity in the $\lambda$-parameter. \hfill\eop

\end{itemize}

\bigskip

\noindent The following easy Lemma gives a bound of each of the three terms appearing in
the right-hand side (\ref{eq:Jalpha}). Since the first term $-\lambda \int_0^t e^{-2^{-j}s} \tilde{V}(\frac{\alpha_s}{\lambda}) \, ds$ is used to compensate possibly large values 
of the other terms, it is not obvious a priori that all these terms are bounded
separately. The idea of the proofs, always based on the same comparison  trick
already used in the proof of Corollary \ref{cor:main} (see below), is interesting in itself, and will be re-used several times
in the next paragraph. Let us first introduce a useful notation.

\begin{Definition}  let $h_{h_0,g}$ be the solution of the KPZ equation (\ref{eq:KPZdowntoj}) with data $h_0,g$. \end{Definition}

\begin{Lemma} \label{lem:ok}
Assume $(h_0, g)$ are bounded data,  and let
$X=X^{\alpha^*}$ be the solution of the SDE (\ref{eq:SDE}), where $h=h_{h_0,g}$. Let also $f\in C([0,T],{\cal W}^{1,\infty})$.

\begin{itemize}
\item[(i)] 
\BEA &&  \esper^{0,x} \left[ \int_0^t e^{-2^{-j}s} f_{t-s}(X_s) ds\right]\lesssim 
\lambda^{-1/(\beta-1)} \left\{ e^{-2^{-j}t} P_-^{-1}\circ P_+\left( |||\lambda^{1/(\beta-1)} h_0|||_{{\cal W}_j^{1,\infty;P_+}}(x)\right)+ \right. \nonumber\\
&& \qquad \left. +
P_-^{-1}\circ P_+\left(|||\lambda^{1/(\beta-1)} (g+f) |||_{\widetilde{\cal W}_j^{1,\infty;P_{\pm}}([0,t])}(x)
\right)+ P_-^{-1}\circ P_+\left( |||\lambda^{1/(\beta-1)}g|||_{\widetilde{\cal W}_j^{1,\infty;P_{\pm}}([0,t])}(x) \right) \right\};  \nonumber\\ \label{eq:bound-g} \EEA
  in particular, letting $f=\pm |g|$, 
\BEA && \esper^{0,x} \left[ \int_0^t e^{-2^{-j}s} |g_{t-s}(X_s)| ds\right] \lesssim
\lambda^{-1/(\beta-1)} \left\{ e^{-2^{-j}t} P_-^{-1}\circ P_+\left( |||\lambda^{1/(\beta-1)} h_0|||_{{\cal W}_j^{1,\infty;P_+}}(x)\right)+ \right. \nonumber\\
&& \qquad \left. +
P_-^{-1}\circ P_+\left(|||2\lambda^{1/(\beta-1)} g |||_{\widetilde{\cal W}_j^{1,\infty;P_{\pm}}([0,t])}(x)
\right) \right\}; \EEA

\item[(ii)]
\BEA &&  \esper^{0,x} \left[e^{-2^{-j}t} |h_0(X_t)| \right] \lesssim \lambda^{-1/(\beta-1)}
\left\{ 
e^{-2^{-j}t} P_-^{-1}\circ P_+\left(  |||2\lambda^{1/(\beta-1)} h_0|||_{{\cal W}_j^{1,\infty;P_+}}(x) \right) \right. \nonumber\\ && \qquad \qquad \left. 
+ P_-^{-1}\circ P_+\left(|||\lambda^{1/(\beta-1)}g|||_{\widetilde{\cal W}_j^{1,\infty;P_{\pm}}([0,t])}(x)\right) \right\};  \label{eq:bound-h0} \EEA

\item[(iii)]

\BEA  && 0\le \esper^{0,x} \left[ \int_0^t e^{-2^{-j}s} \lambda\tilde{V}\left(\frac{\alpha^*_s}{\lambda} \right) ds\right] \lesssim \lambda^{-1/(\beta-1)} \left\{  e^{-2^{-j}t} P_-^{-1}\circ P_+\left( |||2\lambda^{1/(\beta-1)}h_0|||_{{\cal W}_j^{1,\infty;P_+}}(x)\right)+ \right. \nonumber\\ && \qquad\qquad \left.
P_-^{-1}\circ P_+\left( |||2 \lambda^{1/(\beta-1)} g|||_{\widetilde{\cal W}_j^{1,\infty;P_{\pm}}([0,t])}(x)\right) \right\}. \label{eq:bound-V} \EEA
\end{itemize}
\end{Lemma}

{\bf Proof.}
By definition,
\BEA && \esper^{0,x} \left[ \int_0^t e^{-2^{-j}s} f_{t-s}(X_s) ds\right] = 
\esper^{0,x} \left[ \int_0^t e^{-2^{-j}s}\left(-\lambda\tilde{V}\left(\frac{\alpha^*_s}{\lambda}\right)+(f+g)(t-s,X_s)\right)ds \right. \nonumber\\ && \left. \qquad \qquad +e^{-2^{-j}t} h_0(X_t) \right]-h(t,x) 
\le  h_{h_0,f+g}(t,x)-h(t,x)\EEA
which gives the first bound. Varying the initial condition instead, one gets for a general
$\bar{h}_0\in W^{1,\infty}$
\BEA e^{-2^{-j}t} \esper^{0,x} \left[(\bar{h}_0-h_0)(X_t)\right] &=& 
\esper^{0,x} \left[ \int_0^t e^{-2^{-j}s}\left(-\lambda\tilde{V}\left(\frac{\alpha^*_s}{\lambda}\right)+g(t-s,X_s)\right)ds+e^{-2^{-j}t} \bar{h}_0(X_t) \right]-h(t,x) \nonumber\\
&\le &  h_{\bar{h}_0,g}(t,x)-h(t,x). \label{eq:barh0} \EEA
Letting $\bar{h}_0=h_0+|h_0|$  gives an estimate for $e^{-2^{-j}t}\esper^{0,x} [|h_0(X_t)|]$.
Thus one also has a bound for $\esper^{0,x} \left[ \int_0^t e^{-2^{-j}s} \lambda\tilde{V}\left(\frac{\alpha^*_s}{\lambda}\right)ds \right]=-h(t,x)+ \esper^{0,x} \left[ \int_0^t e^{-2^{-j}s} g_{t-s}(X_s) ds\right]+e^{-2^{-j}t} \esper^{0,x} h_0(X_t).$
\hfill \eop

\bigskip


\subsection{Convergence results, existence and unicity of the solution}


Let $(h_0^{(n)})_{n\ge 1}$, resp. $(g^{(n)})_{n\ge 1}$ be sequences of 
bounded data in ${\cal W}^{1,\infty}$, resp. $C([0,T],{\cal W}^{1,\infty})$, converging
locally uniformly to $h_0$, resp. $g$, as in Definition \ref{intro:def:Wj1inftyPsolution}.

\medskip

\noindent Our two main technical Lemmas are the following.

\begin{Lemma} \label{lem:Hnn'}
Let $K\subset\R^d$ be a fixed compact. Then $h^{(n)}\big|_{[0,T]\times K}$ is
a Cauchy sequence of $C([0,T]\times K)$. 
\end{Lemma}

{\bf Proof.} Let $H^{(n,n')}:=h^{(n)}-\mu h^{(n')}$, for $n,n'\in\N^*$, $\mu\in[0,1)$
(ultimately we want to take the limit $\mu\to 1$). Then $H^{(n,n')}(t,x)\le \sup_{\alpha} J_{\alpha}(t,x)$, where
\BEQ J_{\alpha}(t,x)=\esper\left[ \int_0^t e^{-2^{-j}s} \left( -\lambda(1-\mu) \tilde{V}(\frac{\alpha_s}{\lambda})+ (g^{(n)}-\mu g^{(n')})(X_s^{\alpha}) \right)\, ds + e^{-2^{-j}t}
(h_0^{(n)}-\mu h_0^{(n')})(X_t^{\alpha}) \right],\EEQ
$X^{\alpha}$ being as usual the solution of the stochastic differential equation
$dX_s=\alpha_s ds+ dB_s$. We rewrite $g^{(n)}-\mu g^{(n')}$ as $(1-\mu)g^{(n')}+
(g^{(n)}-g^{(n')})$, and similarly, $h_0^{(n)}-\mu h_0^{(n')}$ as $(1-\mu)h_0^{(n')}+
(h_0^{(n)}-h_0^{(n')})$, and split the additive factor $-\int_0^t e^{-2^{-j}s} \lambda(1-\mu)\tilde{V}(\frac{\alpha_s}{\lambda}) ds$ into
two equal parts. We need here {\em polynomial} quasi-norms, and take $d_-=2,d_+>2+\frac{\beta-1}{\beta}d$.  One term is easily bounded, 

\BEA  && (1-\mu) \esper\left[ \int_0^t  e^{-2^{-j}s} \left( -\half\lambda \tilde{V}(\frac{\alpha_s}{\lambda})+  g^{(n')}(X_s^{\alpha}) \right)\, ds + e^{-2^{-j}t}
 h_0^{(n')}(X_t^{\alpha}) \right]  \nonumber\\
 && \qquad\qquad \lesssim \lambda^{-1/(\beta-1)} (1-\mu) \left( |||\lambda^{1/(\beta-1)} g|||^{d_+/2}_{\widetilde{\cal W}_j^{0,\infty;P_{\pm}}([0,t])}(x)+ |||\lambda^{1/(\beta-1)} h_0|||^{d_+/2}_{{\cal W}_j^{0,\infty;P_+}}(x) \right).  \nonumber\\ \label{eq:2.37} \EEA
The other term must still be split into two by using $\esper[\ \cdots\ ]=\esper[ {\bf 1}_{\Omega} \cdots\ ]+ \esper[ {\bf 1}_{\Omega^c} \cdots\ ]$, where $\Omega$ is the
event: $\sup_{0\le s\le t} |X_s-x|\ge L$. First, on $\Omega^c$,
\BEA &&  \esper\left[ {\bf 1}_{\Omega^c} \left\{   \int_0^t e^{-2^{-j}s} \left( -\half \lambda(1-\mu) \tilde{V}(\frac{\alpha_s}{\lambda})+ (g^{(n)}- g^{(n')})(X_s^{\alpha}) \right)\, ds + e^{-2^{-j}t}
(h_0^{(n)}- h_0^{(n')})(X_t^{\alpha}) \right\} \right] \nonumber\\
&& \qquad \qquad \le t ||g^{(n)}-g^{(n')}||_{\infty,[0,t]\times B(x,L)} + ||h_0^{(n)}-h_0^{(n')}||_{\infty,B(x,L)}. \label{eq:2.38} \EEA

Finally one must bound
\BEA  I_{\Omega}(x) & :=&   \esper\left[ {\bf 1}_{\Omega} \left\{   \int_0^t e^{-2^{-j}s} \left( -\half \lambda(1-\mu) \tilde{V}(\frac{\alpha_s}{\lambda})+ (g^{(n)}- g^{(n')})(X_s^{\alpha}) \right)\, ds + e^{-2^{-j}t}
(h_0^{(n)}- h_0^{(n')})(X_t^{\alpha}) \right\} \right] \nonumber\\
&=& (1-\mu) \esper\left[ {\bf 1}_{\Omega} \left\{   \int_0^t e^{-2^{-j}s} \left( -\half \lambda \tilde{V}(\frac{\alpha_s}{\lambda})+ \frac{(g^{(n)}- g^{(n')})(X_s^{\alpha})}{1-\mu} \right)\, ds + e^{-2^{-j}t}
\frac{(h_0^{(n)}- h_0^{(n')})(X_t^{\alpha})}{1-\mu} \right\} \right] \nonumber\\ \EEA
Roughly speaking, $I_{\Omega}(x)$ is small because the characteristic function
${\bf 1}_{\Omega}$ vanishes except: (i) if $\sup_{0\le s\le t} |B_s|\ge \frac{L}{2}$ is large, an event $\Omega'$ of small probability; or (ii) if $\int_0^t |\alpha_s|ds\ge \frac{L}{2}$ is large, implying that $\int_0^t e^{-2^{-j}s} \tilde{V}(\alpha_s) \, ds$ is large (one
retrieves the usual dichotomy). Case (i) is the easier one: we get by Cauchy-Schwarz's
inequality 
\BEA &&(1-\mu) \esper\left[ {\bf 1}_{\Omega\cap \Omega'} \left\{   \int_0^t e^{-2^{-j}s} \left( -\half \lambda \tilde{V}(\frac{\alpha_s}{\lambda})+ \frac{(g^{(n)}- g^{(n')})(X_s^{\alpha})}{1-\mu} \right)\, ds + e^{-2^{-j}t}
\frac{(h_0^{(n)}- h_0^{(n')})(X_t^{\alpha})}{1-\mu} \right\} \right] \nonumber\\
&& \le  (1-\mu) \proba[\Omega']^{1/2} \left( \esper\left[  \left\{   \int_0^t e^{-2^{-j}s} \left( -\half \lambda \tilde{V}(\frac{\alpha_s}{\lambda})+ \frac{(g^{(n)}- g^{(n')})(X_s^{\alpha})}{1-\mu} \right)\, ds + e^{-2^{-j}t}
\frac{(h_0^{(n)}- h_0^{(n')})(X_t^{\alpha})}{1-\mu} \right\}^2 \right] \right)^{1/2}
\label{eq:CS}
\nonumber\\
\EEA
The probability $\proba[\Omega']$ is very small, of order $O(e^{-cL^2/t})$, while
the expectation  in (\ref{eq:CS}) is bounded (up to a coefficient 2 to a certain power)
like the square of the solution $h^{(n,n')}$ of (\ref{eq:KPZdowntoj}) with data $(\frac{h_0^{(n)}-h_0^{(n')}}{1-\mu}, \frac{g^{(n)}-g^{(n')}}{1-\mu})$. All together, using Theorem
\ref{th:main}, we see that (\ref{eq:CS}) is bounded by $O(e^{-cL^2/t})$, times
\BEQ 
 \left( \frac{\lambda^{1/(\beta-1)}}{1-\mu} \right)^{\frac{d_+}{2}-1} \left( |||h_0^{(n)}-h_0^{(n')} |||_{{\cal W}_j^{0,\infty;P_+}}^{d_+/2}(x) + |||g^{(n)}-g^{(n')}|||^{d_+/2}_{\widetilde{\cal W}_j^{0,\infty;P_{\pm}}([0,t])}(x) \right). \label{eq:hnn'} \EEQ 
 
Consider now case (ii). We again split the additive factor $-\half\int_0^t e^{-2^{-j}s} \lambda(1-\mu)\tilde{V}(\frac{\alpha_s}{\lambda}) ds$ into
two equal parts and use the fact, already shown in the proof of Theorem \ref{th:main},  that\\  $-\frac{1}{4}\int_0^t e^{-2^{-j}s} \lambda\tilde{V}(\frac{\alpha_s}{\lambda}) ds \lesssim - \lambda^{-\beta/(\beta-1)} e^{-2^{-j}t} \min\left( \frac{L^2}{t},
\left( \frac{L^{\beta}}{t} \right)^{1/(\beta-1)} \right) $ deterministically, to which
we must add (up to a coefficient 2 to a certain power) the bound (\ref{eq:hnn'}) for $h^{(n,n')}$ found above.

We  shall now choose the parameter $\mu$ and the distance $L\gg 2^{j/2}$ in a near-optimal way to bound $\sup_{[0,t]\times B(x,2^{j/2})} H^{(n,n')}$. Let $\eps>0$. First  we choose $\mu$
so that  $\lambda^{-\beta/(\beta-1)} e^{-2^{-j}t} \min\left( \frac{L^2}{t},
\left( \frac{L^{\beta}}{t} \right)^{1/(\beta-1)} \right)>>$ (\ref{eq:hnn'}), in order
that (ii) does not contribute. This may be arranged uniformly in $n,n'$ for 
\BEQ 1-\mu \approx C(\lambda,t)
 \left( |||h_0 |||_{{\cal W}_j^{0,\infty;P_+}}(x) +
  |||g|||_{\widetilde{\cal W}_j^{0,\infty;P_{\pm}}([0,t])}(x) \right)
L^{-2\beta/(\beta-1)d_+} \label{eq:1-mu}, \EEQ where $C(\lambda,t)$ is a function of $\lambda$ and $t$ (which 
will change from line to line). For such a value of $\mu$, (i) is bounded by\\
\BEQ C(\lambda,t) \left( |||h_0 |||_{{\cal W}_j^{0,\infty;P_+}}(x) + |||g|||_{\widetilde{\cal W}_j^{0,\infty;P_{\pm}}([0,t])}(x) \right) e^{-c'L^2/t} \label{eq:Hnn'(i)}. \EEQ
Then we choose $L$ large enough so that 
\BEA && \max\left( C(\lambda,t) \left( |||h_0 |||_{{\cal W}_j^{0,\infty;
P_+}}(x) + |||g|||_{\widetilde{\cal W}_j^{0,\infty;P_{\pm}}([0,t])}(x) \right)
 e^{-c'L^2/t}, \right. \nonumber\\ && \qquad\qquad \left.  (1-\mu) \left( |||g|||^{d_+/2}_{\widetilde{\cal W}_j^{0,\infty;P_{\pm}}([0,t])}(x)+
  |||h_0|||^{d_+/2}_{{\cal W}_j^{0,\infty;P_+}}(x) \right), (1-\mu)\sup_{[0,t]\times B(x,2^{j/2})} |h^{(n')}| \right) \nonumber\\
&& \lesssim C(\lambda,t) \max\left(  \left( |||h_0 |||_{{\cal W}_j^{0,\infty;
P_+}}(x) + |||g|||_{\widetilde{\cal W}_j^{0,\infty;P_{\pm}}([0,t])}(x) \right)
 e^{-c'L^2/t}, \right. \nonumber\\ && \qquad\qquad \left.  L^{-2\beta/(\beta-1)d_+}  \left( |||g|||^{1+d_+/2}_{\widetilde{\cal W}_j^{0,\infty;P_{\pm}}([0,t])}(x)+
  |||h_0|||^{1+d_+/2}_{{\cal W}_j^{0,\infty;P_+}}(x) \right)  \right) \nonumber\\ \EEA
(taking into account the term coming from (\ref{eq:2.37})) is $<\eps/2$.  
 Finally we choose $n_0$ large enough so that, for all $n,n'\ge n_0$,
$t ||g^{(n)}-g^{(n')}||_{\infty,[0,t]\times B(x,L)} + ||h_0^{(n)}-h_0^{(n')}||_{\infty,B(x,L)}$ (coming from (\ref{eq:2.38})) is also $<\eps/2$.  All together we have proved: $\sup_{[0,t]\times B(x,2^{j/2})} |h^{(n)}-h^{(n')}|\le (1-\mu) \sup_{[0,t]\times B(x,2^{j/2})} |h^{(n')}|+ \sup_{[0,t]\times B(x,2^{j/2})} |H^{(n,n')}| <\eps.$ \hfill\eop

 \bigskip

\begin{Lemma}[non-explosion for random characteristics] \label{lem:non-explosion}
Let $X:=X^{\alpha^*}$ be the path optimizing the Hamilton-Jacobi-Bellman 
problem. Choose $P_{\pm}(z)=z^{d_{\pm}}$ with $d_-=2, d_+>2+\frac{\beta-1}{\beta}d$.
Then 
\BEQ \proba[|X_t-x|\ge L] \lesssim C\left(\lambda,t,|||h_0|||_{{\cal W}_j^{0,\infty;P_+}}(x), |||g|||_{\widetilde{\cal W}_j^{0,\infty;P_{\pm}}([0,t])}(x) \right) L^{-2\beta/(\beta-1)d_+}. \EEQ
\end{Lemma}

{\bf Proof.} We let $\bar{h}_0(x):=h_0(x)+{\bf 1}_{|x|\ge L}$, 
remark that $\proba[|X_t-x|\ge L]=\esper^{0,x}\left[ (\bar{h}_0-h_0)(X_t)\right]$, and  use
(\ref{eq:barh0}), getting $e^{-2^{-j}t} \proba[|X_t-x|\ge L]\le h_{\bar{h}_0,g}(t,x)-h(t,x)$.
Next we rewrite $h_{\bar{h}_0,g}(t,x)-h(t,x)$ as $H(t,x)-(1-\mu)h(t,x)$, with
$H(t,x):=h_{\bar{h}_0,g}(t,x)-\mu h(t,x)$ ($0\le \mu<1$), and proceeds to bound the function $H$ in a similar way as in the previous Lemma: letting $\bar{X}:=\bar{X}^{\bar{\alpha}^*}$ be the
path optimizing the Hamilton-Jacobi-Bellman problem with data $(\bar{h}_0,g)$,  and
rewriting $\bar{h}_0-\mu h_0$ as $(1-\mu)h_0+(\bar{h}_0-h_0)$, we find
\BEQ H(t,x) \le  (1-\mu) (H_1(t,x)+H_2(t,x)),\EEQ
where
\BEQ H_1(t,x):=\esper\left[ \int_0^t  e^{-2^{-j}s} \left(-\half
\lambda \tilde{V}(\frac{\bar{\alpha}^*_s}{\lambda})+ g(\bar{X}_s) \right)\, ds + e^{-2^{-j}t}
h_0(\bar{X}_t) \right] \EEQ 
and
\BEQ H_2(t,x):= \esper\left[ -\int_0^t  e^{-2^{-j}s} \half
\lambda  \tilde{V}(\frac{\bar{\alpha}^*_s}{\lambda})\, ds + \frac{1}{1-\mu} e^{-2^{-j}t}
(\bar{h}_0-h_0)(\bar{X}_t) \right]. \label{eq:H2} \EEQ

By Theorem \ref{th:main}, the function $(1-\mu)H_1$ is bounded by a constant times\\
 $(1-\mu)  (\lambda^{1/(\beta-1)})^{\frac{d_+}{2}-1} \ \left\{  |||g|||^{d_+/2}_{\widetilde{\cal W}_j^{0,\infty;P_{\pm}}([0,t])}(x)+ e^{-2^{-j}t}
  |||h_0|||^{d_+/2}_{{\cal W}_j^{0,\infty;P_+}}(x) \right\}$. So let us find an upper
  bound for the function $(1-\mu) H_2$.  Since $\bar{h}_0-h_0$ vanishes on $B(x,L)$, we may insert the characteristic
  function ${\bf 1}_{\Omega}$ inside the expectation, $\esper[\, \cdot\,]\to
  \esper[{\bf 1}_{\Omega} \, \cdot\,]$ in (\ref{eq:H2}), where $\Omega$ is
  the event: $|\bar{X}_t-x|\ge L$. Comparing with the proof of Lemma \ref{lem:Hnn'}, see in
  particular (\ref{eq:hnn'}),  we see
  that $(1-\mu)H_2(t,x)$ is bounded by 
 \BEA &&  O(e^{-cL^2/t})\left( \frac{\lambda^{1/(\beta-1)}}{1-\mu} \right)^{\frac{d_+}{2}-1} \left( |||h_0 |||_{{\cal W}_j^{0,\infty;P_+}}^{d_+/2}(x) + |||g|||^{d_+/2}_{\widetilde{\cal W}_j^{0,\infty;P_{\pm}}([0,t])}(x) \right)  \nonumber\\
  && \qquad +
   \left\{  \left( \frac{\lambda^{1/(\beta-1)}}{1-\mu} \right)^{\frac{d_+}{2}-1} \left( |||h_0 |||_{{\cal W}_j^{0,\infty;P_+}}^{d_+/2}(x) + |||g|||^{d_+/2}_{\widetilde{\cal W}_j^{0,\infty;P_{\pm}}([0,t])}(x) \right) - c\lambda^{-\beta/(\beta-1)} e^{-2^{-j}t} \min\left( \frac{L^2}{t},
\left( \frac{L^{\beta}}{t} \right)^{1/(\beta-1)} \right) \right\}.
\label{eq:non-explosion} \nonumber\\ \EEA
 We choose $\mu$ in the same way as in the proof of Lemma \ref{lem:Hnn'}, see (\ref{eq:1-mu}), in such a way that the  second term in (\ref{eq:non-explosion}) does not
 contribute. As a result,  
 $(1-\mu)H_2(t,x)$ is bounded by\\ $C(\lambda,t) \left( |||h_0 |||_{{\cal W}_j^{0,\infty;P_+}}(x) + |||g|||_{{\cal W}_j^{0,\infty;P_{\pm}}([0,t])}(x) \right) e^{-c'L^2/t}$ as
 in (\ref{eq:Hnn'(i)}). Finally, we get
 \BEQ \proba[|X_t-x|\ge L]\le H(t,x)+(1-\mu)|h(t,x)| \le C(\lambda,t) \left( |||h_0 |||^{1+d_+/2}_{{\cal W}_j^{0,\infty;P_+}}(x) +
  |||g|||^{1+d_+/2}_{{\cal W}_j^{0,\infty;P_{\pm}}([0,t])}(x) \right)
L^{-2\beta/(\beta-1)d_+}. \EEQ  \hfill \eop

\bigskip
From Lemma \ref{lem:Hnn'} we know that all sequences $(h^{(n)})$ converge uniformly on any compact to the same function $h$. By the stability principle for continuous viscosity solutions 
(see e.g. \cite{Bar}), the limit is a solution of the KPZ equation with data
$(h_0,g)$. Furthermore, from Corollary \ref{cor:main}, we know that the functions
$h^{(n)}, n\ge 1$ are locally uniformly differentiable, i.e. for all compact $K\subset\R^d$, and all $T>0$, $\sup_{n\ge 1}\sup_{t\le T} \sup_K |\nabla h_t^{(n)}|
<C(K,T)$. We deduce using Lemma 3.11 of \cite{Unt-KPZ1} that $h$ is classical.
This concludes the proof of Theorem 1 and Theorem 2.


\section{Appendix}



\subsection{Integration lemmas}


\begin{Lemma} \label{lem:beta}
Let $\Phi:\R^d\to\R_+$ be a rotationally invariant, positive, smooth function such
that $\Phi$ is a decreasing function of the module of its argument,  and $\int \Phi(x)\, dx=1$. Then
\BEQ \int dy\,  \Phi(x-y)|f(y)|dy\le C f^*(x)\EEQ
for some universal constant $C$ (depending only on $d$).
\end{Lemma}

{\em Proof.} By abuse of notation, we write $\Phi(r)$ for the value of $\Phi$ on
the sphere $S_r=\partial B(0,r)$; note that $\Phi'(r)\le 0$. Then
\BEA \int dy\,  \Phi(x-y)|f(y)|dy &=& \int dr\,  \Phi(r) \int_{S_r} dy |f(y)|=
\int  dr\,  |\Phi'(r)| \int_{B(x,r)} dy \, |f(y)| \nonumber\\
&\le & f^{\sharp}(x) \int dr\,  |\Phi'(r)| \int_{B(x,r)} dy=f^{\sharp}(x) \int dy\,  \Phi(x-y)
=f^{\sharp}(x),\EEA
where $f^{\sharp}(x):=\sup_{r>0} \fint_{B(x,r)} |f|$ is the supremum of the local
averages of $|f|$ around $x$ (see Introduction or \cite{Unt-KPZ1}, section 3.1). We conclude
using \cite{Unt-KPZ1}, Lemma 3.2. \hfill\eop

\medskip

The above lemma applies to the {\em generalized heat kernels with exponent $\beta$}, 
\BEQ \Phi_t^{\beta}(x-y):=c \frac{e^{-(|x-y|^{\beta}/t)^{1/(\beta-1)}}}{t^{d/\beta}},
\label{eq:Phibeta} \EEQ
where $c^{-1}=c^{-1}(\beta):=\int \Phi_1^{\beta}(y) dy$ is a normalization constant (in
particular, $\Phi^2$ is the usaul heat kernel). 
Note also the following elementary truncated integral estimates,
\BEQ \int_{|y|>A} dy\,  \Phi(y) |f(x-y)|\le f^{\sharp}(x)\int_{|y|>A} dy\, \Phi(y).\EEQ In particular, an
integration by parts yields for $A\gg t^{1/\beta}$
\BEA  \int_{|y|>A} dy\,  \Phi_t^{\beta}(y) |f(x-y)| &\lesssim&  f^{\sharp}(x) A^{d-1} (t/A)^{1/(\beta-1)}
\Phi_t^{\beta}(A)\nonumber\\
&\lesssim& f^{\sharp}(x) e^{-\half (A^{\beta}/t)^{1/(\beta-1)}}.\EEA

\medskip

The next Lemma is a generalization of Lemma \ref{lem:integration0}. Instead of
covering $\R^d$ by cubes with side $2^{j/2}$, it is convenient to choose slightly
{\em smaller} cubes, with side $c2^{j/2}$ ($0<c<1$), and let $\widetilde{\D}^j:=\left\{\tilde{\Del}=[c(k_1-\half) 2^{j/2}, c(k_1+\half)2^{j/2}]\times\cdots\times
[c(k_d-\half)2^{j/2}, c(k_d+\half)2^{j/2}], \ k_1,\ldots,k_d\in\Z\right\}$,  chosen in such a way that the origin cube $\tilde{\Del}_0:=[-\half c2^{j/2},\half c2^{j/2}]\times\cdots\times [-\half c2^{j/2},\half c2^{j/2}]$ is entirely included in the ball $B(0,2^{j/2})$.
With this definition we have:

\begin{Lemma} \label{lem:local-integration}
Let $f\in L^{\infty}_{loc}(\R^d)$ such that $\locsup^j(f)^*(x)<\infty$; $X^x$ an
$\R^d$-valued random variable depending on $x\in\R^d$;  and $w:\R^d\to\R_+$ a positive weight
function such that, for all $x\in\R^d$, and $\tilde{\Del}\in\widetilde{\D}^j$,
\BEQ \left(\sup_{y\in\tilde{\Del}} w(y) \right) \ \cdot\  \proba[X^x-x\in\tilde{\Del}]\le C\Phi(\min_{y\in\tilde{\Del}}|y|) \EEQ
for some kernel $\Phi$  as in Lemma \ref{lem:beta}. Then, for every $\tau\ge 0$,
\BEQ e^{\tau\Del} \left( x\mapsto \esper[ {\bf 1}_{[X^x-x|>2^{j/2}} w(X^x-x) |f(X^x)|
] \right)(x) \lesssim C(2^{-j})^{d/2} \locsup^j(f)^*(x).\EEQ
\end{Lemma}

{\bf Proof.} One finds

\BEA e^{\tau\Del} \left( x\mapsto\esper[ {\bf 1}_{[X^x-x|>2^{j/2}} w(X^x-x) |f(X^x)|] \right) &\le& e^{\tau\Del} \left(x\mapsto   \sum_{\tilde{\Del}\not=\tilde{\Del}_0}
\left(\sup_{y\in\tilde{\Del}} w(y)\right) \ \cdot\  \proba[X^x-x\in\tilde{\Del}] \sup_{\tilde{\Del}} |f(x+\cdot)| \right) \nonumber\\
&\le & C\sum_{\tilde{\Del}\not=\tilde{\Del}_0} \Phi(\min_{y\in\tilde{\Del}}|y|)  e^{\tau\Del} (x\mapsto \sup_{\tilde{\Del}}|f(x+\cdot)|) \nonumber\\
&\lesssim & C(2^{-j})^{d/2} \int dz\, \Phi(c|z|) e^{\tau\Del}(\locsup^j(f))(x+z) \nonumber\\
&\lesssim &  C(2^{-j})^{d/2} \locsup^j(f)^*(x)\EEA
where $c>0$ is some constant.

\hfill \eop


\subsection{Estimates in pointwise quasi-norms}


We sketch in this paragraph how to adapt our main results (Theorem 1 and Theorem 2
in the Introduction) to the case when the initial condition $h_0$ and the right-hand side
$g$ have possibly large fluctuations on regions with small relative volume, so that e.g. 
$|||\locsup^j h_0|||_{{\cal H}^P}(x)$, $|||2^{j/2}\locsup^j|\nabla h_0|\, |||_{{\cal H}^{P_{\pm}}}(x)$ are very large or infinite. Typically, one may imagine that e.g. 
$h_0$ is constructed in the following way (see Example in section 3.1 of \cite{Unt-KPZ1}). Define $P(h_0)$ to be identically equal to $c_1>0$ outside some union of annuli $\cup_{k\ge 0} B_k$, where $B_k:=B(0,2^k+2^{k\gamma})\setminus B(0,2^k)$ for some $\gamma\in(-\infty,1)$, and $P(h_0)\big|_{B_k}(x):=
c_1+(c_2 (2^k)^{d(1-\gamma')}-c_1) \chi(\frac{|x|-2^k}{2^{k\gamma}})$
($\gamma'\ge \gamma$), where $\chi:[0,1]\to[0,1]$ is some smooth 'bump' function such
that $\chi\big|_{[0,\frac{1}{4}]\cup [\frac{3}{4},1]}\equiv 0$, $\chi\big|_{[\frac{1}{3},\frac{2}{3}]} \equiv 1$  (one may actually choose exponents $\gamma'_k\ge\gamma_k$
depending on $k$). If $\gamma'<0$, then $\locsup^j(P(h_0))(x)\approx
|x|^{d(1-\gamma')}>>|x|^d$ in the annuli $(B_k)_{k\ge 0}$,  in contradiction
with the growth condition (\ref{eq:growth-condition}),  so that $|||\locsup^j h_0|||_{{\cal H}^P}(x)=\infty$; while, on the other hand, $(P(h_0))^*(x)=O(c_1)+O(c_2)$ is
finite, as is easy to check. Even under such circumstances,  it is possible to prove existence and
uniqueness of the solution, and prove good estimates, in much larger functional spaces defined
by   {\em pointwise quasi-norms} instead of {\em local quasi-norms}. 

In principle, the whole idea boils down to this: {\em pointwise} quasi-norms, denoted
by the sub-index "point", $|||\, \cdot \, |||_{{\cal H}^P_{point}}(x)$, $|||\, \cdot\, 
|||_{{\cal H}^{P_{\pm}}([0,t])}(x)$, $|||\, \cdot\, |||_{{\cal W}_{j,point}^{1,\infty;P}}(x)$, $|||\, \cdot\, |||_{{\cal W}_{j,point}^{1,\infty;P_{\pm}}([0,t])}(x)$ (see \S 0.3) are
derived from {\em local} quasi-norms by moving the suprema out of the averaging
"star" operation. Some of this work has already been done in some details in
\cite{Unt-KPZ1} (see end of \S 3.4):

\begin{itemize}
\item[(i)]  replace $\locsup^j(f)^*=(\locsup^j(f))^*$ by $\locsup^j (f^*)$. 
Thus $|||\locsup^j h_0|||_{{\cal H}^P}(x)$ in (\ref{eq:Wj1inftyP}) shoud be replaced
by $\locsup^j\left(|||h_0|||_{{\cal H}^P}\right)(x)$.
\item[(ii)] rewrite $2^{j/2} \locsup^j|\nabla f|(x)$ as $\sup_{\vec{\eps},\vec{\eps}'\in B(0,1)} |\del_{\vec{\eps},\vec{\eps}'} f(x)|$, where $\del_{\vec{\eps},\vec{\eps}'}f(x):=
\frac{f(x+2^{j/2}\vec{\eps})- f(x+2^{j/2}\vec{\eps}')}{|\vec{\eps}-\vec{\eps}'|}$.
Thus $|||2^{j/2}\locsup^j|\nabla h_0|\, |||_{{\cal H}^P}(x)$ in (\ref{eq:Wj1inftyP})
shoud be replaced by $\sup_{\vec{\eps},\vec{\eps}'} |||\del_{\vec{\eps},\vec{\eps'}} h_0
|||_{{\cal H}^P}(x)$.
\end{itemize}

We must still deal with the time-dependent norms, see (\ref{eq:Wj1inftyP0t}).
From the time-discretization used in the proof of Theorem \ref{th:main}, it is
clear (iii) that the integral $\int_0^t e^{-2^{-j}s} P_-^{-1}\circ P_+\left(
|||2^j \Locsup^j g(t-s,\cdot)|||_{{\cal H}^{P_+}}(x)\right)$ appearing in the formula for $|||\Locsup^j g|||_{{\cal H}^{P_{\pm}}([0,t])}(x)$ (see  (\ref{intro:def:g-norm}) and (\ref{eq:Wj1inftyP0t})) should
be replaced e.g. with $\locsup^j \left( \sup_{n\ge 1} \frac{t}{n} \sum_{k=0}^{n-1} e^{-2^{-j}kt/n}
P_-^{-1}\circ P_+\left( |||2^j g(t-\frac{kt}{n},\cdot)|||_{{\cal H}^{P_+}(x)} \right)\right)$. 
Finally, combining (ii) and (iii), the integral$\int_0^t e^{-2^{-j}s} P_-^{-1}\circ P_+\left(
|||2^{3j/2} \Locsup^j |\nabla g|(t-s,\cdot)|||_{{\cal H}^{P_+}(x)}\right)$ appearing in the formula for $|||2^{j/2} \Locsup^j |\nabla g|\ |||_{{\cal H}^{P_{\pm}}([0,t])}(x)$
should be replaced with $\sup_{n\ge 1} \sup_{\vec{\eps},\vec{\eps}'\in B(0,1)}
\frac{t}{n} \sum_{k=0}^{n-1} e^{-2^{-j}kt/n}
P_-^{-1}\circ P_+\left( |||2^{3j/2} \del_{\vec{\eps},\vec{\eps}'} g(t-\frac{kt}{n},\cdot)|||_{{\cal H}^{P_+}}(x) \right)$.

\bigskip

Apart from the fact that the very definition of the pointwise quasi-norms appears to
be quite lengthy, this substitution makes the actual proof of estimates like
those appearing in Theorem 2 much more technical. The main reason is that
we cannot use Lemma \ref{lem:local-integration} any more. Instead, we have the much more
difficult to deal with

\begin{Lemma} \label{lem:pointwise-integration}
Let $X$ be an $\R^d$-valued random variable,  $r^{d-1}d\sigma_r$ the surface measure on the sphere
$S_r=\{|x|=r\}\subset\R^d$ of radius $r$, and $u:\R^d\to\R$ a $d$-times continuously differentiable function. Then
\BEQ \esper[u(X)]\lesssim  \sup_{y\in B(0,2^{j/2})} u(y)+ \sum_{k=0}^d (2^{-j/2})^{d-k}
\int_{2^{j/2}}^{+\infty} dr\  \proba[|X|>r] \ \int_{S_r}
|\nabla^k u(x)| r^{d-1}d\sigma_r(x). \label{eq:pointwise-integration} \EEQ
\end{Lemma}

{\bf Proof.}  We use a finite partition of unity, $1\equiv \chi_0+\chi_1+\chi_2$ such that: $\chi_0,\chi_1,\chi_2\ge 0$; $\supp(\chi_0)\subset B(0,2^{j/2})$, $\supp(\chi_i)= \{rx, x\in \Omega_i; r\ge 2^{j/2-1}\}$, $i=1,2$
where $\Omega_1$, resp. $\Omega_2$ are closed subsets of the unit sphere containing  the northern, resp. southern hemisphere.
We choose two smooth diffeomorphisms (e.g. spherical coordinates up to normalization) $\phi_i: \Omega_i\to [0,1]^{d-1}$, $i=1,2$, let $\Phi_i:\supp(\chi_i)\to \R_+\times [0,1]^{d-1}, 
rx \mapsto (r,\phi_i(x))$ be its "suspension"  with inverse $\Phi_i^{-1}(r,y)=r\phi_i^{-1}(y)$,
and assume  $\tilde{\chi}_i:=\chi_i\circ\Phi_i^{-1}$ satisfy the natural hypotheses $||\partial_r^k \nabla_y^{l} \tilde{\chi}_i||_{\infty}=
O((2^{-j/2})^k)$ for arbitrary $k\ge 0$ and multi-index $l$. Finally, we may assume by a density argument that
$X$ has a smooth density $f$, and let $\tilde{f}_i:=|J|f\circ\Phi_i^{-1}$ be the  densities of the transfered
variables $\tilde{X}_i:=X\circ\Phi_i^{-1}$,  where
 $J(r,y):=\left|\frac{d\Phi_i^{-1}(r,y)}{d(r,y)}\right|$ is  the Jacobian; similarly, we let $\tilde{u}_i:=u\circ\Phi_i^{-1}$.

We first write $\esper[u(X)]\le \sup_{B(0,2^{j/2})} u + \sum_{i=1,2}\esper[u(X)\chi_i(X)]$ and
\BEQ \esper[u(X)\chi_i(X)]=\int \chi_i(x)u(x)f(x)dx=\int \tilde{\chi_i}(r,y) \tilde{u}_i(r,y)\tilde{f}_i(r,y) dr dy.\EEQ
Inspired by the one-dimensional integration-by-parts formula, 
\BEQ \esper[u(Y)]=u(0)+\int_0^{+\infty} u'(y)\proba[Y>y]
dy \qquad (Y\ge 0) \label{eq:IPP}, \EEQ
 we integrate by parts with respect to each of the $d$ coordinates in $\R_+\times[0,1]^{d-1}$. Since there are no
boundary terms, one gets (letting $\tilde{X}_i^1,\ldots,\tilde{X}_i^d$ be the coordinates of $\tilde{X}_i$ and
$y=(y_2,\ldots,y_d)$)
\BEQ \int \tilde{\chi_i}(r,y) \tilde{u}_i(r,y)\tilde{f}_i(r,y) dr dy=\int \left[ \prod_{m=2}^{d} \left(\frac{1}{r}\partial_{y_m}\right)\partial_r
(\tilde{\chi_i} \tilde{u}_i) \right] (r,y) \proba[\tilde{X}_i^1>r,\tilde{X}_i^2>y_2,\ldots, \tilde{X}_i^d>y_d]\  r^{d-1}dr dy.\EEQ
Now $\proba[\tilde{X}_i^1>r,\tilde{X}_i^2>y_2,\ldots, \tilde{X}_i^d>y_d]\le \proba[|X|>r]$; 
the derivative with respect to $r$ produces a factor $O(2^{-j/2})$, resp. $O(1)$ when applied to the cut-off $\tilde{\chi}_i$,
resp. to the function $\tilde{u}_i$,
while normalized angular derivatives $\frac{1}{r}\partial_{y_m}$ yield factors $O(r^{-1})\lesssim 2^{-j/2}$, resp. $O(1)$. All together one gets the result.
\hfill \eop

\bigskip

\noindent If we want to apply this Lemma instead of Lemma \ref{lem:local-integration} in
the course of the proof of Theorem \ref{th:main}, then we must consider data $h_0,g$ such that
not only $h_0\in {\cal W}_{j,point}^{1,\infty;P_+}$, resp. $g\in \widetilde{\cal W}_{j,point}^{1,\infty;P_{\pm}}([0,T])$, but also their space derivatives, $2^{j|\kappa|/2}\nabla^{\kappa} h_0$, resp. $2^{j|\kappa|/2}\nabla^{\kappa}g$, to all orders $|\kappa|\le d$ (mind that the natural scaling factors $2^{j|\kappa|/2}$, read off (\ref{eq:pointwise-integration}), are mainly decorative in the polynomial case, but not in the exponential
case!) This defines  new functional spaces ${\cal W}^{d+1,\infty;P}_{j,point}$, ${\cal W}^{d+1,\infty;P_{\pm}}_{j,point}$.   

\medskip

\noindent The conclusion of this discourse is the following: assuming
(i) (exponential case) $P_+(z)=e^{\lambda z}$, $P_-(z)=e^{c\lambda z}$ with
$c\in(0,1)$ small enough; (ii) (polynomial case) $P_{\pm}(z)=z^{d_{\pm}}$ with
$d_-\ge 1$, $d_+-d_-$ large enough; and $d'$ (the time-exponent used in the 
passage from time-dependent ${\cal W}$-quasi norms to the modified $\widetilde{\cal W}$-quasi norm) large enough,  one should be able to prove that
\BEQ |||h_t|||_{{\cal W}^{1,\infty;P_-}_{j,point}}(x) \lesssim e^{-2^{-j}t}
|||h_0|||_{{\cal W}_{j,point}^{d+1,\infty;P_-}}(x) + |||g|||_{\widetilde{\cal W}_{j,point}^{d+1,\infty;P_{\pm}}([0,t]]}(x). \label{eq:fin} \EEQ
 
\noindent Finally, Schauder estimates (e.g. in the form proved in \cite{Unt-Bur1}) make it
possible to show pointwise bounds for $|\nabla^{\kappa} h(t,x)|$, $|\kappa|=2,\ldots,d+1$ in terms of the local suprema of $h$, $\nabla h$ and the space-derivatives of 
the data $u_0,g$ up to order $d+1$, which probably implies a bound of the type of (\ref{eq:fin}) for $|||h_t|||_{{\cal W}_{j,point}^{d+1,\infty;P_-}}(x)$ with a 
larger loss of regularity. 
 


\begin{thebibliography}{99}






\bibitem{AmoBen} L. Amour, M. Ben-Artzi. {\em Global existence and decay for viscous Hamilton-Jacobi equations},
Nonlinear Analysis, Theory, Methods and Applications {\bf 31}, 621--628 (1998).

\bibitem{Bar} G. Barles. {\em Solutions de viscosit\'e et \'equations elliptiques du deuxi\`eme ordre},
(graduate course given at the University of Tours, available at: www.lmpt.univ-tours.fr/\~\,barles/Toulcours.pdf (1997).

\bibitem{BenBenLau} S. Benachour, M. Ben-Artzi, P. Lauren\c cot. {\em Sharp decay estimates and vanishing viscosity
for diffusive Hamilton-Jacobi equations}, Adv. Differential Equations {\bf 14}, 1--25 (2009).

\bibitem{BenKarLau} S. Benachour, G. Karch, P. Lauren\c cot. {\em Asymptotic profiles of solutions to viscous
Hamilton-Jacobi equations}, J. Math. Pures Appl. {\bf 83}, 1275--1308 (2004).

\bibitem{BenLau} S. Benachour, P. Lauren\c cot. {\em Global solutions to viscous Hamilton-Jacobi equations with
irregular initial data}, Comm. Partial Diff. Eq. {\bf 24}, 1999-2021 (1999).



\bibitem{Bre} H. Br\'ezis. {\em Op\'erateurs maximaux monotones et semi-groupes de contractions dans les
espaces de Hilbert}, North-Holland (1973).

\bibitem{BricGawKup} J. Bricmont, K. Gawedzki, A. Kupiainen. {\em KAM theorem and quantum field theory},
 Comm. Math. Phys. {\bf 201}, 699-727 (1999).

\bibitem{ComYos} F. Comets, N. Yoshida. {\em Directed polymers in random environment are diffusive at weak
disorder}, Ann. Prob. {\bf 34}, 1746--1770 (2006).

\bibitem{Doo} J. L. Doob. {\em Stochastic processes}, John Wiley and Sons (1967).


\bibitem{Eva} L. C. Evans. {\em Partial differential equations}, Graduate Studies in Mathematics {\bf 19}, AMS (199).

\bibitem{Fle} W. H. Fleming, R. W. Rishel. {\em Deterministic and stochastic optimal control}, Applications
of mathematics, Stochastic modelling and applied probability, vol. 1, Springer (1975).

\bibitem{Frie}
A. Friedman. {\em Partial differential equations of parabolic type}, Prentice-Hall(1964).


\bibitem{GilGueKer} B. H. Gilding, M. Guedda, R. Kersner. {\em The Cauchy problem for $u_t=\Del u+|\nabla u|^q$},
J. Math. Anal. Appl. {\bf 284}, 733-755 (2003).




\bibitem{Ito} K. Ito. {\em Existence of solutions to the Hamilton-Jacobi-Bellman equation under quadratic
growth conditions}, J. Diff. Eq. {\bf 176}, 11-28 (2001).

\bibitem{Joh} K. Johansson, {\em Shape fluctuations and random matrices}, Comm. Math. Phys. {\bf 209},
437--476 (2000).



\bibitem{KS} I. Karatzas, S. Shreve. {\it Brownian motion and stochastic calculus},
Springer-Verlag (1991).

\bibitem{Kry} N. V. Krylov. {\em Controlled diffusion processes}, Applications of mathematics, vol. 14, Springer
(1980).

\bibitem{LauSou} P. Lauren\c cot, P. Souplet. {\em On the growth of mass for a viscous Hamilton-Jacobi
equation}, Jour. Anal. Math. {\bf 89}, 367--383 (2003).

\bibitem{Lie} G. M. Lieberman. {\em Second order parabolic differential equations},
World Scientific (1996).


\bibitem{RY} D. Revuz, M. Yor.  {\it Continuous martingales and Brownian motion}, Springer (1999).

\bibitem{Rud} W. Rudin. {\em Real and complex analysis}, McGraw-Hill, New York (1966).

\bibitem{Sme} I. Smears. {\em Hamilton-Jacobi-Bellman equations. Analysis and numerical analysis}, research
report available on www.math.dur.ac.uk/Ug/projects/highlights/PR4/Smears\_HJB\_report.pdf


\bibitem{Trie} H. Triebel. {\em Theory of function spaces II}, Monographs in 
Mathematics {\bf 84}, Birkh\"auser (1992).


\bibitem{Unt-KPZ1}  J. Unterberger. {\it PDE estimates for multi-dimensional
KPZ equation}, arXiv:1307.1980.

\bibitem{Unt-Bur1} J. Unterberger. {\it Global existence and smoothness for solutions
of viscous Burgers equation. (1) The bounded case}, arXiv:1503.05145.

\bibitem{Unt-Bur2} J. Unterberger. {\it Global existence and smoothness for solutions
of viscous Burgers equation. (2) The unbounded case: a characteristic flow
study},  arXiv:1510.01539.

\bibitem{XJW} X.-J. Wang. {\em Schauder estimates for elliptic and parabolic equations}, Chin. Ann. Math.
 {\bf 27B(6)}, 637-642 (2006).


\end{thebibliography}
\end{document}